\documentclass[a4paper,12pt,leqno,twoside]{article}
\usepackage[english]{babel} 
\usepackage{amsmath}
\usepackage{amsfonts}
\usepackage{amssymb}
\usepackage{latexsym}
\usepackage{xypic}
\usepackage[all]{xy}
\usepackage{stmaryrd} 
\usepackage{theorem}  
\usepackage[mathscr]{eucal}
\usepackage{eufrak}
\usepackage{epsf}
\usepackage[dvips]{graphicx}

\usepackage{epsfig}

\usepackage{latexsym}  
\usepackage[usenames]{color}

\usepackage{fancyhdr}
\usepackage{mathrsfs} 

\usepackage{helvet}

\renewcommand{\phi}{\varphi}
\renewcommand{\theta}{\vartheta}
\renewcommand{\rho}{\varrho}
\newcommand{\ep}{\epsilon}
\newcommand{\etaep}{{\eta_\epsilon}}
\newcommand{\lra}{\longrightarrow}

\newcommand{\tpi}{{\tilde\pi}}
\renewcommand{\bigr}[1]{{\big(#1\big)}}
\renewcommand{\Bigr}[1]{{\Big(#1\Big)}}
\renewcommand{\biggr}[1]{{\bigg(#1\bigg)}}
\newcommand{\bigc}[1]{{\big\{#1\big\}}}

\newcommand{\com}{\mathbb{C}}
\newcommand{\rea}{\mathbb{R}}

\newcommand{\integer}{\mathbb{Z}}
\newcommand{\integergz}{\mathbb{Z}_{>0}}
\newcommand{\nat}{\mathbb{N}}

\newcommand{\reagz}{\mathbb{R}_{>0}}

\newcommand{\dist}[2]{\mathrm{dist}(#1,#2)}
\newcommand{\dpart}[1]{\delta_{\partial#1}}

\newcommand{\de}[2]{\delta_{#1}(#2)}
\newcommand{\norm}[3]{|\!|#1|\!|^{#2}_{#3}}

\newcommand{\gsect}[3]{{S_{#1,#2,#3}}}
\newcommand{\csect}[3]{{S_{#1\pm #2,#3}}}

\newcommand{\ball}[2]{B(#1,#2)}
\newcommand{\band}[2]{B_{#1}^{#2}}

\newcommand{\Op}{\mathrm{Op}}

\newcommand{\opxsac}{\mathrm{Op}^c_{\xsa}}

\newcommand{\cov}[1]{\mathrm{Cov}_{sa}(#1)}

\newcommand{\xsa}{{X_{sa}}}

\newcommand{\xsalf}{{X_{sa,l\!f}}}

\newcommand{\rsa}{{\rea^2_{sa}}}

\newcommand{\Ho}[3][]{\mathcal{H}\mathrm{om}_{#1}(#2,#3)}

\newcommand{\RH}[3][]{\mathit{R}\mathcal{H}\mathit{om}_{#1}(#2,#3)}
\newcommand{\Rh}[3][]{\mathit{RHom}_{#1}(#2,#3)}
\newcommand{\M}{\mathcal{M}}
\newcommand{\ot}{\mathcal{O}^t}
\newcommand{\otxsa}{\mathcal{O}^t_\xsa}

\newcommand{\otwsa}{\mathcal{O}^t_{W_{sa}}}
\renewcommand{\O}{\mathcal{O}}
\newcommand{\dbt}{\mathcal{D}b^t}

\newcommand{\dbxr}{\mathcal{D}b_{X_\rea}}

\newcommand{\dbtxsa}{\mathcal{D}b^t_{\xsa}}

\newcommand{\mcC}{\mathcal{C}}

\newcommand{\D}{\mathcal{D}}
\newcommand{\DX}{{\D_X}}
\newcommand{\dbdx}[1][]{D^b_{#1}(\DX)}
\newcommand{\bdc}[2][]{D^b_{#1}(#2)}
\renewcommand{\mod}{\mathrm{Mod}}
\renewcommand{\M}{\mathcal{M}}

\renewcommand{\O}{\mathcal{O}}

\newcommand{\ch}{\mathrm{char}}
\newcommand{\dmod}[1]{\mathrm{Mod}(\mathcal D_{#1})}
\newcommand{\cdmod}[1]{\mathrm{Mod}_{coh}(\mathcal D_{#1})} 
\newcommand{\hdmod}[1]{\mathrm{Mod}_{h}(\mathcal D_{#1})} 
\newcommand{\lM}{\mathcal{M}(*a)}
\newcommand{\sol}{\mathscr{S}ol}
\newcommand{\solt}{\mathscr{S}ol^t}
\newcommand{\cF}{\mathscr F}

\newtheorem{thm}{Theorem}[subsection]
\newtheorem{df}[thm]{Def\mbox{}inition}
\newtheorem{cor}[thm]{Corollary}
\newtheorem{prop}[thm]{Proposition}
\newtheorem{lem}[thm]{Lemma}

\newtheorem{rem}[thm]{Remark}

\newtheorem{notation}[thm]{Notation}

\numberwithin{equation}{section}

\newcommand{\proofend}{\hfill $\Box$ \vspace{\baselineskip}\newline}
\title{\bf{An existence theorem for tempered solutions of
    $\mathcal{D}$-modules on complex curves}}
\date{\today}

\author{Giovanni Morando}

\newpage

\begin{document}

\maketitle

\thispagestyle{empty}





















\begin{abstract}
Let $X$ be a complex curve, $\xsa$ the subanalytic site associated to
$X$, $\M$ a holonomic $\D_X$-module. Let
$\otxsa$ be the sheaf on $\xsa$ of tempered holomorphic functions
and $\sol(\M)$ (resp. $\solt(\M)$) the complex of holomorphic
(resp. tempered holomorphic) solutions of $\M$. We prove that the
natural morphism
$$ H^1(\solt(\M))\lra H^1(\sol(\M))$$
is an isomorphism. As a consequence, we prove that $\solt(\M)$ is
$\rea$-constructible in the sense of sheaves on $\xsa$. Such a result
is conjectured by M. Kashiwara and P. Schapira in \cite{arx} in any
dimension.

\vspace{4mm}


\end{abstract}





\tableofcontents


\addcontentsline{toc}{section}{\textbf{Introduction}}
\section*{Introduction}
 \nopagebreak 

The problem of existence for ordinary linear differential equations (and even
non-linear) is classical and the
litterature presents many results on this subject. In particular,
existence theorems for solutions with growth conditions have been
obtained by many authors such as J.-P. Ramis and Y. Sibuya
(\cite{ramis_sibuya}), B. Malgrange (\cite{malgr_birk}) and N. Honda
(\cite{honda}). In \cite{ramis_sibuya} and
\cite{malgr_birk}, the authors proved existence for functions with
Gevrey-type growth conditions at the origin on sectors of sufficiently
small amplitude. Using similar techniques, in \cite{honda}, the author
proved existence for ultra-distributions with support on $\rea_{\geq0}$.

The functional spaces considered in \cite{ramis_sibuya} and
\cite{malgr_birk} correspond to sheaves on the real blow-up at the origin
of $\com$. Essentially they are sheaves on the unit circle. Indeed,
growth conditions did not allow a global sheaf theoretical approach.

Nonetheless tempered distributions were a basic tool in M. Kashiwara's
functorial proof of the Riemann-Hilbert correspondence in \cite{kashi79} and
\cite{rims}. In order to use tempered distributions functorially,
M. Kashiwara introduced the new functor $T\mathcal H om$ of tempered
cohomology. Such a functor represented the
first step in a different approach to sheaves which, through
\cite{moderate_and_formal}, led to the full use of sheaves on sites in
\cite{ast}. Indeed in \cite{ast}, M. Kashiwara and P. Schapira combined
classical analytical results of S. \L ojasiewicz (\cite{loj_studia},
see also \cite{malgr_ideals}) with sheaves on sites. They
realized tempered distributions, tempered $\mathcal C^\infty$
functions and Whitney $\mathcal C^\infty$ functions as sheaves on the
subanalytic site. They also defined tempered holomorphic functions
$\otxsa$ as the complex of the solutions of the Cauchy-Riemann system
in the space of tempered distributions. 
 
In a subsequent paper, \cite{arx}, M. Kashiwara and P. Schapira extended the
notion of microsupport of sheaves to the subanalytic site. In this way
they established the framework for the study of tempered holomorphic solutions of
$\D$-modules. They also gave an example which is the starting point of
the study of tempered holomorphic solutions of an irregular ordinary
differential equation.

Given a complex analytic manifold $X$, we denote by $\Op_\xsa^c$ the
category of relatively compact subanalytic open subsets of $X$ and by
$\xsa$ the subanalytic site, that is the site whose underlying
category is $\Op_\xsa^c$ and whose coverings are the finite
coverings. We denote by $\mod(k_X)$ (resp. $\mod(k_\xsa)$) the
category of sheaves of $k$-modules on the site $X$ (resp. $\xsa$). Let
$\rho:X\lra\xsa$ be the natural morphism of sites.

Given a $\D_X$-module $\M$, it is natural to compare 
$$ \sol\M   :=  R\rho_*\mathrm R\mathcal Hom_{\D_X}\big(\M,\O_X\big) $$ 
and 
$$ \solt\M  :=  \mathrm R\mathcal Hom_{\rho_!\D_X}\big(\rho_!\M,\ot_\xsa\big)  \  $$  
(for the definition of $\rho_!$, see Section \ref{temp_hol}).

Along his proof of the Riemann-Hilbert correspondence, M. Kashiwara
proved that, if $\M$ is a regular $\D_X$-module, then
$\solt\M\overset{\sim}{\lra}\sol\M$.

In \cite{arx}, the authors studied $\solt\M$ comparing it to $\sol\M$,
for a particular example on a complex curve $X$. 

In the present paper, we go into the study of $\solt(\M)$ for $\M$ a
holonomic $\D$-module on a complex curve $X$. In particular we prove 
an existence theorem for tempered solutions of ordinary differential
equations in the subanalytic topology, thus refining the classical
results on small open sectors. Such a result has two consequences.

First, we obtain that the natural morphism
\begin{equation}\label{eq_H1_intro} 
H^1(\solt\M)\lra H^1(\sol\M)  
\end{equation}
is an isomorphism.

Second, we prove that the complex $\solt(\M)$ is $\rea$-constructible
in the sense of \cite{arx}. In that paper the authors conjectured such
a result in any dimension.

Our results being on a complex curve, it is natural to look for
extensions of them in higher dimensions. 
In \cite{sabbah_ast}, C. Sabbah conjectured and widely
developed the higher dimensional version of Hukuhara-Turrittin's
Theorem. Recently Y. Andr{\'e} announced the proof of Sabbah's
conjecture. Such results would be at the base of a possible extensions
of our results.

The contents of the present paper are subdivided as follows.

In \textbf{Section \ref{SEC_AN_GEOM}}, we briefly review subanalytic sets recalling the
classical results that we will need. We study in detail relatively
compact subanalytic open subsets of $\rea^2$. 
We give a decomposition of
$U\in\Op_\rsa^c$ using sets biholomorphic to open sectors. Such a
result will be essential in Section 3.

In \textbf{Section \ref{SEC_TEMP_HOL}}, we recall definitions and
basic results of sheaves on the subanalytic site and
tempered holomorphic functions on a complex curve. In Subsection 2.2 we prove a
result concerning the composition of a tempered holomorphic function
and a biholomorphism. 

In \textbf{Section \ref{SEC_EXIST_THM}} we consider an open disc $X\subset\com$
centered at $0$ and
\begin{equation}\label{eq_diff_op_intro} 
 P:=z^N\frac d{dz}I_m+A(z) \  ,
\end{equation}
where $m\in\integer_{>0}$, $N\in\nat$, $ A\in gl(m,\O_\com(X))$ and
$I_m$ is the identity matrix of order $m$. The aim of this section is
to study the solvability of $P$ in the space of tempered holomorphic
functions on $U\in\opxsac$ with $0\in\partial U$. We prove that there
exist an open neighborhood $W\subset\com$ of $0$ and a finite
subanalytic covering $\{U_j\}_{j\in J}$ of $U\cap W$
such that for any $g_j\in\otxsa(U_j)^m$ there exists
$u_j\in\otxsa(U_j)^m$ such that $Pu_j=g_j$, for any $j\in J$. We start
the section by recalling Hukuhara-Turrittin's
Theorem which is a basic tool in the study of ordinary differential
equations.


In \textbf{Section \ref{SEC_DMOD}}, we deal with $\D_X$-modules on a complex analytic curve
$X$. We begin by recalling some classical results on
$\D_X$-modules. In Subsection \ref{SUBSEC_H1}, we prove a first
consequence of the results of Section \ref{SEC_EXIST_THM}, that is,
\eqref{eq_H1_intro} is an isomorphism. In Subsection \ref{SUBSEC_R-C}
we prove a second consequence of the results of Section
\ref{SEC_EXIST_THM}, that is, $\solt(\M)$ is $\rea$-constructible in
the sense of sheaves on $\xsa$.

We thank P. Schapira for proposing this problem to our attention and
for many fruitful discussions and A. D'Agnolo for many useful remarks.

\section{Subanalytic sets}  \label{SEC_AN_GEOM}
In the f\mbox{}irst subsection, we recall the definition and some
classical results on subanalytic sets. In the second subsection we
focus on relatively compact subanalytic open subsets of $\rea^2$. We
prove some results mixing the complex and the real analytic structure
on $\rea^2$. Indeed, we describe the local structure of relatively
compact subanalytic open subsets of $\rea^2$ via biholomorphic images
of open sectors (Theorem \ref{thm_cov_suban}). 

\subsection{Review on subanalytic sets}

Let $M$ be a real analytic manifold, $\mathcal A$ the sheaf of
real-valued real analytic functions on $M$.

\begin{df}
Let $X\subset M$. 
\begin{enumerate}
\item $X$ is said to be \emph{semi-analytic at $x\in M$} if 
 there exists an open neighborhood $W$ of $x$ such that $X\cap
 W=\cup_{i\in I}\cap_{j\in J}X_{i,j}$ where $I$ and $J$ are f\mbox{}inite sets
 and either $X_{i,j}=\{y\in W;\ f_{i,j}(y)>0\}$ or $X_{i,j}=\{y\in W;\
 f_{i,j}(y)=0\}$ for some $f_{i,j}\in\mathcal A(W)$. 
 $X$ is said \emph{semi-analytic} if it is semi-analytic at each $x\in
M$.
\item $X$ is said \emph{subanalytic} if for any $x\in M$ there exist an open
  neighborhood $W$ of $x$, a real analytic manifold $N$ and a
  relatively compact semi-analytic set $A\subset M\times N$ such that
  $\pi(A)=X\cap W$, where $\pi: M\times N\to M$ is the projection. 
\item Let $N$ be a real analytic manifold. A map $f:X\to N$ is
  said \emph{subanalytic} if its graph, 
$$\Gamma_f:=\bigc{(x,y)\in X\times N;y=f(x)} \ , $$ 
is subanalytic in $M\times N$. 
\end{enumerate}\end{df}

Given $X\subset M$, denote by $\overset{\circ}{X}$ (resp. $\overline X$, $\partial X$) the interior (resp. the closure, the boundary).

\begin{prop}[See \cite{bier_mil}]
Let $X$ and $Y$ be subanalytic subset of $M$. Then $X\cup Y$, $X\cap
Y$, $\overline X$, $\overset{\circ}{X}$ and $X\setminus Y$ are
subanalytic. Moreover the connected components of $X$ are subanalytic,
the family of connected components of $X$ is locally f\mbox{}inite and $X$ is
locally connected.
\end{prop}

Def\mbox{}inition \ref{def_ccd}, Theorem \ref{thm_existence_ccd} and
Proposition \ref{prop_limit_subanalytic_function} below are stated and
proved in \cite{coste} for the more general case of o-minimal structures.

\begin{df}[Cylindrical Cell Decomposition]\label{def_ccd}
Let $n\in\integergz$. A \emph{cylindrical cell decomposition} (\emph{ccd} for short) $\bigc{C_k}_{k\in K}$ of
$\rea^n$ is a f\mbox{}inite partition of $\rea^n$ into
subanalytic sets $C_k$ obtained inductively on $n$ in the following way. The sets $C_k$
are called \emph{cells}.
\begin{itemize}
\item [$n=1$:] The cells def\mbox{}ining a ccd of $\rea$ are open intervals
  $]a,b[$ or points $\{c\}$, where $a\in\rea\cup\{-\infty\}$,
  $b\in\rea\cup\{+\infty\}$, $a<b$, and $c\in\rea$.
\item [$n>1$:] A ccd
  $\bigc{D_h}_{h\in H}$ of $\rea^n$ is given by a ccd $\bigc{C_k}_{k\in K}$ of $\rea^{n-1}$, 
  $l_k\in\nat$ and subanalytic continuous functions
  $$\zeta_{k,1},\ldots,\zeta_{k,l_k}:C_k\to\rea$$ 
  such that, for any $x\in C_k$,
  $\zeta_{k,j}(x)<\zeta_{k,j+1}(x)$, $j=1,\ldots,l_k-1$ $(k\in K)$.

The cells $D_h$ are the graphs of $\zeta_{k,j}$,
$$\phantom{\quad\quad 1\leq j\leq l_k} \Gamma_{\zeta_{k,j}}:=\big\{\bigr{x,\zeta_{k,j}(x)}\in C_k\times\rea\big\} \quad\quad (1\leq j\leq l_k) \ , $$
and the sets
\begin{equation}\label{eq_band} 
\big\{(x,y)\in C_k\times\rea;\zeta_{k,j}(x)<y<\zeta_{k,j+1}(x)\big\} 
\end{equation}
for $0\leq j\leq l_k$, where $\zeta_{k,0}=-\infty$ and $\zeta_{k,l_k+1}=+\infty$.
\end{itemize}
\end{df}

\begin{thm}[See \cite{coste}, Theorem 2.10]\label{thm_existence_ccd}
\sloppy Let $A_1,\ldots,A_d$ be relatively compact subanalytic subsets of
$\rea^n$. There exists a cylindrical cell decomposition of $\rea^n$
adapted to each $A_j$. That is, each $A_j$ is a union of cells.
\end{thm}

\begin{prop}[See \cite{coste}, Theorem 3.4]\label{prop_limit_subanalytic_function}
Let $Z$ be a subanalytic subset of $\rea^n$. The following properties
are equivalent.
\begin{enumerate}
\item $Z$ is closed and bounded.
\item Every subanalytic continuous map $\zeta:]0,1[\to Z$ extends by
  continuity to a map $[0,1[\to Z$. 
\item For any subanalytic continuous function $\zeta:Z\to\rea$, $\zeta(Z)$ is
  closed and bounded.
\end{enumerate}
\end{prop}

For Theorem \ref{thm_loj_ineq} below, see \cite[Theorem 6.4]{bier_mil}.

\begin{thm}[\L ojasiewicz's Inequality.]\label{thm_loj_ineq}
Let $M$ be a real analytic ma\-ni\-fold, $K\subset M$. Let
$f,g:K\to\rea$ be subanalytic functions with compact graphs. If
$f^{-1}\bigr{\{0\}}\subset g^{-1}\bigr{\{0\}}$, then there exist $c,r\in\reagz$
such that, for any $x\in K$,
$$ |f(x)|\geq c|g(x)|^r \ .$$ 
\end{thm}

For Theorem \ref{thm_csl} below, see \cite[Proposition 8.2.3]{som}.

\begin{thm}[Curve Selection Lemma.]\label{thm_csl}
Let $Z$ be a subanalytic subset of $M$ and let $z_0\in\overline Z$. Then there exists an ana\-ly\-tic map
$$ \gamma:]-1,1[\longrightarrow M \ ,$$
such that $\gamma(0)=z_0$ and $\gamma(t)\in Z$ for $t\neq0$.
\end{thm}

\subsection{Subanalytic subsets of $\rea^2$}\label{SUBSECTION_COV_SUBAN}

\begin{notation}
\sloppy Given a real analytic manifold $M$, we denote by $\Op_{M_{sa}}$
(resp. $\Op^c_{M_{sa}}$) the category of subanalytic open (resp. relatively compact subanalytic open) subsets of $M$. 
\end{notation}

Let 
$$\begin{array}{rrcl}
\tpi:  &  \rea_{\geq0}\times]-\pi,3\pi[  &  \lra  &  \rea^2 \\
  & (\rho,\theta)  &  \longmapsto  &  \rho e^{i\theta}  \ .
\end{array}   $$

One has that, given $U\in\Op_\rsa^c$ with $0\notin U$, $\tpi^{-1}(U)$ is a
subanalytic open subset of $\reagz\times]-\pi,3\pi[$, relatively
compact in $\rea^2$.

For $R\in\reagz$, $\eta,\xi:[0,R]\lra]-\pi,3\pi[$ subanalytic
continuous functions such that $\eta(\rho)<\xi(\rho)$, for any
$\rho\in]0,R[$, denote
$$ \band\eta\xi:=\Big\{(\rho,\theta)\in\,]0,R[\,\times\,]-\pi,3\pi[\,;\,\eta(\rho)<\theta<\xi(\rho)\Big\} \ . $$

Remark that $\overline{\band\eta\xi}\subset[0,R]\times]-\pi,3\pi[$.

\begin{prop}\label{prop_blowup}
Let $U\in\Op_\rsa^c$, $0\in\partial U$. 
There exists an open neighborhood $W\subset\rea^2$
of $0$, such that $U\cap W$ is a finite union of sets of the form
$\tpi(\band\eta\xi)\cap W$.
\end{prop}

\emph{Proof.}
The set $\tpi^{-1}(U)$ is a subanalytic open subset of
$\reagz\times]-\pi,3\pi[$, relatively compact in $\rea^2$. Let
$\ep\in\reagz$, $\ep<\pi$. Take a cylindrical cell decomposition of
$\rea^2$ adapted to
$$\tpi^{-1}(U)\cap\Bigr{\reagz\times\big]-\ep,2\pi+\ep\big[} \ . $$ 
The conclusion follows.
\proofend

For $z\in\com$ and $\epsilon\in\reagz$, denote by $B(z,\epsilon)$
the open ball of center $z$ and radius $\epsilon$. 

Let us introduce semi-analytic arcs and prove an easy result which
states the local equivalence between semi-analytic arcs and graphs
of subanalytic functions.

\begin{df}
Let $\gamma:]-1,1[\longrightarrow\rea^2$ be an analytic map,
$\delta\in\reagz$ such that $\gamma|_{[0,\delta]}$ is
injective. We call
$$\Gamma:=\gamma\bigr{]0,\delta[}$$ 
a \emph{semi-analytic arc with an endpoint at $\gamma(0)$}.
\end{df}

Recall that, given a function $\eta$, we denote by $\Gamma_\eta$ the graph of $\eta$.

\begin{lem}\label{lem_graph=arc}
Let $R\in\reagz$, $\eta:[0,R[\to\rea$ a subanalytic continuous
map. There exist $\delta\in\reagz$ and an analytic map
$\gamma:]-1,1[\longrightarrow\rea^2$ such that
$\gamma(0)=\big(0,\eta(0)\big)$ and 
\begin{equation}\label{eq_graph=arc} 
\gamma\big(]-1,1[\setminus\{0\}\big)=\Gamma_\eta\cap\big(]0,\delta[\times\rea\big)
\  .
\end{equation}
In particular, there exist a semi-analytic arc $\Gamma$ with an endpoint at
$\bigr{0,\eta(0)}$ and an open neighborhood $W\subset\rea^2$ of $\big(0,\eta(0)\big)$, such that 
$$  \Gamma_\eta\cap W=\overline\Gamma\cap W \ .$$
\end{lem}

\emph{Proof.} 
Let $p_1:\rea^2\to\rea$ be the projection on the first
coordinate.

By Theorem \ref{thm_csl} there exists an analytic map
$\gamma:]-1,1[\lra\rea^2$ such that $\gamma(0)=\bigr{0,\eta(0)}$ and
\begin{equation}\label{eq_gamma_subset_Gamma}
\gamma\Bigr{]-1,1[\setminus\{0\}}\ \subset\ \Gamma_\eta\!\setminus\!\big\{\bigr{0,\eta(0)}\big\} \ .
\end{equation} 

Remark that we can suppose that $\gamma|_{[0,1[}$ and $\gamma|_{]-1,0]}$ are injective. 
Since $\gamma\bigr{]-1,1[}$ is arcwise-connected,
$p_1\Bigr{\gamma\bigr{]-1,1[}}$ is arcwise-connected as well. Hence, since
$\{0\}\subsetneqq p_1\Bigr{\gamma\bigr{]-1,1[}}\subset\rea_{\geq0}$, there exists $\delta\in\reagz$
such that $p_1\bigr{\gamma\bigr{]-1,1[}}=[0,\delta[ $.

Further, by \eqref{eq_gamma_subset_Gamma}, 
\begin{equation}\label{eq_0notin_p1}
p_1\Bigr{\gamma\bigr{]-1,1[\,\setminus\{0\}}}=]0,\delta[ \ . 
\end{equation}

Let us prove that if $0<x<\delta$, then
$\bigr{x,\eta(x)}\in\gamma\bigr{]-1,1[\setminus\{0\}}$, this will conclude the proof. Let $x\in]0,\delta[$. By \eqref{eq_0notin_p1}, there exists $y\in\rea$ such that
$(x,y)\in\gamma\bigr{]-1,1[\setminus\{0\}}$. By \eqref{eq_gamma_subset_Gamma}, it
follows that $y=\eta(x)$. Hence $\bigr{x,\eta(x)}\in\gamma\bigr{]-1,1[\setminus\{0\}}$.
\nopagebreak
\proofend

Roughly speaking, from Lemma \ref{lem_graph=arc} and Proposition \ref{prop_blowup},
it follows that $(\partial U)\cap W$ is a finite collection of semi-analytic
arcs with an endpoint at $0$.

Let us now introduce biholomorphic images of open sectors. We start
with a well known result on the local nature of holomorphic functions
on $\com$. For Proposition \ref{prop_local_nature_holomorphic} below,
see \cite[Theorem 2.1]{forster}.

\begin{prop}\label{prop_local_nature_holomorphic}
Let $U\subset\com$ be an open neighborhood of $0$,
$\phi:U\to\com$ a non constant holomorphic map such that $0$ is a
zero of order $n$ for $\phi$. 
There exist an open neighborhood $U'\subset U$ of $0$,
$\epsilon\in\reagz$, and a holomorphic isomorphism
$\psi:U'\to\ball0\epsilon$ 
such that, for any $z\in U'$,
$$ \phi|_{U'}(z)=\big(\psi(z)\big)^n \ .$$
\end{prop}

\begin{df}\label{df_sector} 
Let $\alpha,\beta\in\rea$, $r\in\reagz$, $\alpha<\beta$. 
The set
 $$ \gsect\alpha\beta r:=\big\{\rho e^{i\theta}\in\com^\times;\ 0<\rho<r,\ \theta\in]\alpha,\beta[\big\} $$  
is called an {\em open sector of amplitude $\beta-\alpha$ and radius
 $r$} or simply {\em an open sector}.
\end{df}

We will need to stress on the amplitude and the direction of a
  sector so we will also use the following slightly dif\mbox{}ferent
  notation
$$ \csect\tau\eta r  :=  \gsect{\tau-\eta}{\tau+\eta} r  $$
for $\tau\in\rea$ and $\eta,r\in\reagz$. 

\begin{cor}\label{cor_bastardo}
\sloppy Let $U\subset\com$ be an open neighborhood of $0$,
$\phi:U\to\com$ a non constant holomorphic map such that
$\phi(0)=0$. 
\begin{enumerate}
\item There exist $r,\tau\in\reagz$ such that $\overline{\ball0r}\subset U$ and, for any $\theta\in\rea$, 
  $\phi|_{\overline{\csect\theta\tau r}}$ is an injective map.
\item Suppose that, given $\alpha,\beta\in\rea$, there exist
  $\mu,\delta,R\in\reagz$ such that 
$$ \phi\big(]0,\delta[\times\{0\} \big)\subset\gsect{\alpha+\mu}{\beta-\mu}R \ .$$ 
Then, there exist $\eta,r'\in\reagz$ such that 
$$  \phi\big(\csect0\eta{r'}\big)\subset \gsect\alpha\beta R \ .$$
\end{enumerate}
\end{cor}

\emph{Proof.} It is based on Proposition \ref{prop_local_nature_holomorphic} and the
fact that holomorphic isomorphisms are conformal maps.

\proofend

We are now ready to state and prove the main result of this
section. Denote by $\O_\com$ the sheaf of holomorphic functions on $\com$.

\begin{thm}\label{thm_cov_suban} 
Let $U\in\Op^c_\rsa$, $0\in\partial U$. There exist an open
neighborhood $W\subset\com$ of $0$, a f\mbox{}inite
set $J$, open sectors $S_{j,k}$,
$\phi_{j,k}\in\O_\com\big(\overline{S_{j,k}}\big)$ ($j\in J,k=1,2$)
such that 
\begin{enumerate}
\item $\phi_{j,k}(0)=0$ and $\phi_{j,k}|_{\overline{S_{j,k}}}$ is injective $(j\in
  J,k=1,2)$,
\item $$U\cap W=\bigcup_{j\in J}\Big(\phi_{j,1}\big(S_{j,1}\big)\cap\phi_{j,2}\big(S_{j,2}\big)\Big) \ .$$
\end{enumerate}
\end{thm}

\emph{Proof of Theorem \ref{thm_cov_suban}.}
By Proposition \ref{prop_blowup}, it is sufficient to prove the
statement for $U=\tpi(\band\eta\xi)\cap W$, for $W\subset\com$ an open
neighborhood of $0$.

First we need two technical lemmas.

\begin{lem}\label{lem_preliminar1}
Let $S$ be an open sector, $\phi\in\O_\com\bigr{\overline S}$ such that
$\phi(0)=0$ and $\phi|_{\overline S}$ is injective. Suppose that
there exists $\ep\in\reagz$ such that $\phi(S)\cap\ball0\ep$ is
contained in an open sector of amplitude strictly smaller than
$2\pi$.

Then there exist $r\in\reagz$, an open neighborhood
$V\subset\com$ of $0$ and $\zeta_1,\zeta_2:[0,\ep]\to]-\pi,3\pi[$ subanalytic
continuous functions such that, for any $\rho\in[0,\ep]$, $\zeta_1(\rho)<\zeta_2(\rho)$ and
$$ \phi\bigr{S\cap\ball0r}=\tpi\bigr{\band{\zeta_1}{\zeta_2}}\cap
V  \ .$$

\end{lem}

\emph{Proof.} We limit to give a sketch of the proof which is essentially of topological nature.

There exist $\eta\in[0,2\pi]$, $\mu\in\reagz$, $\mu<\pi$, such that $\phi(S)\cap\ball0\ep\subset\csect\eta\mu\ep$.

Remark that $[\eta-\mu,\eta+\mu]\subset\,]\eta-\pi,\eta+\pi[\,\subset\,]-\pi,3\pi[$. Take a ccd of
$\rea^2$ adapted to 
$$ \tpi^{-1}\Bigr{\phi\bigr{S}\cap\ball0\ep}\cap\Bigr{\reagz\times]\eta-\mu,\eta+\mu[ } \ . $$ 
Since, for any $\delta\in\reagz$, $\phi\bigr{S}\cap\ball0{\delta}$ has just one connected component having $0$ in its
boundary, the conclusion follows.
\proofend

\begin{lem}\label{lem_preliminar2}
Let $R\in\reagz$, $\eta:[0,R]\lra]-\pi,3\pi[$ a subanalytic
con\-ti\-nuous map. There exist an open neighborhood $V\subset\com$ of
$0$, $\tau,r,\ep\in\reagz$, $\phi\in\O_\com\bigr{\overline{\ball0r}}$ and
$\zeta_1,\zeta_2:[0,\ep]\to]-\pi,3\pi[$ subanalytic continuous
functions satisfying the following conditions.
\begin{enumerate}
\item $\phi|_{\overline{\gsect{-\tau}{\tau}{r}}}$ is
  injective.
\item For any $\rho\in[0,\ep]$, $-\pi<\zeta_1(\rho)<\eta(\rho)<\zeta_2(\rho)<3\pi$
  and
\begin{equation}\label{eq_phi_upper}
\phi\bigr{\gsect{0}{\tau}{r}}=\tpi\Big(\band\eta{\zeta_2}\Big)\cap V \ ,
\end{equation}
\begin{equation}\label{eq_phi_lower}
\phi\bigr{\gsect{-\tau}0{r}}=\tpi\Big(\band{\zeta_1}\eta\Big)\cap V \ .
\end{equation}
\end{enumerate}
\end{lem}

\emph{Proof.} Remark that it is sufficient to prove the
statement for $\eta|_{[0,\ep']}$ for some $\ep'\in\reagz$, $\ep'<R$. We
set for short $\eta_{\ep'}:=\eta|_{[0,\ep']}$.

\sloppy Since  $\eta(0)\in]-\pi,3\pi[$, there
exist $\ep,\mu_1\in\reagz$, $\mu_1<\pi$, such that
$[\etaep(0)-\mu_1,\etaep(0)+\mu_1]\subset]-\pi,3\pi[$ and
\begin{equation}\label{eq_phi_graph_in_sector}
\Gamma_\etaep\setminus\bigr{0,\etaep(0)} \subset\,
]0,R[\,\times\,]\etaep(0)-\mu_1,\etaep(0)+\mu_1[  \ .
\end{equation}

First, let us show that there exist an open neighborhood $W\subset\com$ of $]-1,1[$, $\phi\in\O_\com(W)$ and $\delta\in\reagz$ such that $\phi(0)=0$ and 
\begin{equation}\label{eq_phi_line}
\phi\Big(\bigr{]-1,1[\,\setminus\{0\}}\times\{0\}\Big)=\tpi\Big(\Gamma_\etaep\cap\big(]0,\delta[\times\rea\big)\Big)  \ .
\end{equation}

By Lemma \ref{lem_graph=arc}, there exist $\delta\in\reagz$ and an analytic map
$\gamma:]-1,1[\to\rea^2$ such that
$\gamma(0)=\bigr{0,\eta(0)}$ and
$$\gamma\bigr{]-1,1[\setminus\{0\}}=\Gamma_\etaep\cap\big(]0,\delta[\times\rea\big)
\ . $$

Since $\tpi\circ\gamma$ is an analytic map, there exist a complex
neighborhood $W$of $]-1,1[$ and $\phi\in\O_\com(W)$ such that
$\phi|_{]-1,1[\times\{0\}}=\tpi\circ\gamma|_{]-1,1[}$. In particular,
$\phi(0)=0$ and
\begin{equation}\label{eq_third}
\phi\Bigr{\bigr{]-1,1[\,\setminus\{0\}}\times\{0\}}=\tpi\Big(\Gamma_\etaep\cap\big(]0,\delta[\times\rea\big)\Big) \ .
\end{equation}
Hence, \eqref{eq_phi_line} follows.

Now, remark that \eqref{eq_phi_graph_in_sector} implies that
\begin{equation}\label{eq_forth}
\tpi\bigr{\Gamma_\etaep\setminus\bigr{0,\etaep(0)}}\subset\csect{\etaep(0)}{\mu_1}{R} \ .
\end{equation} 
Combining \eqref{eq_third} and \eqref{eq_forth}, we have that
\begin{eqnarray}\label{eq_fifth}
\phi\bigr{]0,1[\times\{0\}}  &  \subset  &
\tpi\Bigr{\Gamma_\etaep\cap\bigr{]0,\delta[\times\rea}}  \notag \\
  &  \subset  &  \tpi\bigr{\Gamma_\etaep\setminus\bigr{0,\etaep(0)}}  \notag \\
  &  \subset  &  \csect{\etaep(0)}{\mu_1}R \  .
\end{eqnarray}
Since $[\etaep(0)-\mu_1,\etaep(0)+\mu_1]\subset\,]-\pi,3\pi[$ and
$\mu_1<\pi$, there exists $\mu_2\in\reagz$ such that $\mu_1<\mu_2<\pi$ and
\begin{equation}\label{eq_sixth}
[\etaep(0)-\mu_2,\etaep(0)+\mu_2]\subset]-\pi,3\pi[ \ .
\end{equation} 
Let $r\in\reagz$ be such that $\overline{\ball0r}\subset W$. Then, Corollary
\ref{cor_bastardo} \emph{(ii)} applies and there exist $\tau\in\reagz$ such that, up to shrinking $r$, 
\begin{equation}\label{eq_seventh}
\phi\bigr{\gsect0{\tau}{r}}\subset\csect{\etaep(0)}{\mu_2}R \ .
\end{equation}
Further, by Corollary \ref{cor_bastardo} \emph{(i)}, up to shrinking $\tau$ and $r$, we have that
$\phi|_{\overline{\gsect{0}{\tau}{r}}}$ is injective.

Then, Lemma \ref{lem_preliminar1} applies and there exist $r'\in\reagz$, an open neighborhood
$V\subset\com$ of $0$ and $\zeta_1,\zeta_2:[0,\ep]\to]-\pi,3\pi[$ subanalytic
continuous functions such that, for any $\rho\in[0,\ep]$,
$\zeta_1(\rho)<\zeta_2(\rho)$ and
$$ \phi(\gsect{0}{\tau}{r'})=\tpi(\band{\zeta_1}{\zeta_2})\cap V  \ .$$

Then, \eqref{eq_sixth} and \eqref{eq_seventh} imply that one can chose
$\zeta_1=\etaep$. Hence \eqref{eq_phi_upper} follows.

Clearly, \eqref{eq_phi_lower} can be proved using the same arguments.

\proofend

\emph{End of the Proof of Theorem \ref{thm_cov_suban}.}

As said above, by Proposition \ref{prop_blowup}, it is sufficient to prove the
statement for $U=\tpi(\band\eta\xi)\cap W$, for $W\subset\com$ an open
neighborhood of $0$.

Consider $\band\eta\xi$. 

\sloppy By Lemma \ref{lem_preliminar2}, there exist $\zeta_1,\zeta_2:[0,\ep]\to]-\pi,3\pi[$,
$r,\tau\in\reagz$,
$\phi_1,\phi_2\in\O_\com\big(\overline{\ball0r}\big)$, $V_1,V_2\subset\com$
open neighborhoods of $0$ such that, for any $\rho\in[0,\ep]$, $\eta(\rho)<\zeta_2(\rho)<3\pi$,
$-\pi<\zeta_1(\rho)<\xi(\rho)$, $\phi_1|_{\overline{\gsect{0}{\tau}{r}}}$
$\phi_2|_{\overline{\gsect{-\tau}{0}{r}}}$ are injective and
$$ \tpi\bigr{\band\eta{\zeta_2}}\cap V_1 = \phi_1(\gsect0{\tau}{r}) \ , $$
$$ \tpi\bigr{\band{\zeta_1}\xi}\cap V_2 = \phi_2(\gsect{-\tau}0{r})  \ . $$

We distinguish two cases: $\xi(0)=\eta(0)$ and $\eta(0)<\xi(0)$.

\emph{(i)} Suppose $\xi(0)=\eta(0)$. 

We have that
$$-\pi<\zeta_1(0)<\eta(0)=\xi(0)<\zeta_2(0)<3\pi \ . $$ 
It follows that there exists $\ep'\in\reagz$ such that, for any $\rho\in[0,\ep']$,
$$\zeta_1(\rho)<\eta(\rho)\leq\xi(\rho)<\zeta_2(\rho) \ . $$ 
Hence, considering $\eta,\xi,\zeta_1,\zeta_2$ as restricted to $[0,\ep']$, we have that
$$\band\eta\xi=\band\eta{\zeta_1}\cap\band{\zeta_2}\xi \ . $$

Now, up to take smaller $\tau,\ep'$, we can suppose that
$\tpi(\band{\zeta_1}\xi)$ and $\tpi(\band\eta{\zeta_2})$ are contained
in an open sector of amplitude strictly smaller than $2\pi$. In
particular, $\tpi$ is a bijection on $\band\eta{\zeta_1}\cup\band{\zeta_2}\xi$.
It follows that
$$
\tpi\big(\band\eta\xi\big)=\tpi\big(\band\eta{\zeta_1}\big)\cap\tpi\big(\band{\zeta_2}\xi\big)
\ .$$

Taking $V:=V_1\cap V_2$, the conclusion follows.

\emph{(ii)} Suppose $\eta(0)<\xi(0)$. 

Up to take smaller $\tau$, there
exist $\ep'\in\reagz$ and $\alpha,\beta:[0,\ep']\to\rea$ constant
functions such that, for any $\rho\in[0,\ep']$,
$$
\eta(\rho)<\alpha(\rho)<\zeta_2(\rho)<\zeta_1(\rho)<\beta(\rho)<\xi(\rho)
\ . $$

It follows that, considering $\eta,\xi,\zeta_1,\zeta_2$ as restricted to $[0,\ep']$, 
$$
\band\eta\xi=\band\eta{\zeta_2}\cup\band\alpha\beta\cup\band{\zeta_1}\xi
\ . $$
The conclusion follows.
\nopagebreak
\proofend

Detailing the proof of Theorem \ref{thm_cov_suban}, one can give a
more precise statement in the following way. 

\begin{rem}\label{rem_refined}
Let $U$, $W$, $\phi_{j,k}$ and $S_{j,k}$ as given in Theorem \ref{thm_cov_suban}. Given
$r,\eta\in\reagz$, there exist an open neighborhood $W'\subset W$ of
the origin, a finite set $J'$ and open sectors $S'_{j',k}\subset
S_{j,k}$ ($j'\in J'$) such that the amplitude (resp. the radius) of $S'_{j',k}$ is
smaller than $\eta$ (resp. $r$) and 
$$ U\cap W'=\bigcup_{j\in
  J}\Big(\phi_{j,1}\big(S'_{j,1}\big)\cap\phi_{j,2}\big(S'_{j,2}\big)\Big) \ .$$ 
\end{rem}

\section{Tempered holomorphic functions}\label{temp_hol} \label{SEC_TEMP_HOL}

In the first subsection we recall the definition and some classical
results on the subanalytic site $\xsa$ underlying a complex curve $X$
and sheaves on $\xsa$. In the second subsection
we recall the definition of the subanalytic sheaf of tempered
holomorphic functions. In the third section we prove a result on the
pull back of tempered holomorphic functions through biholomorphims.
We refer to \cite{arx} and \cite{ast} for the first and the second subsection.

Throughout this section, $X$ will be a complex analytic curve.

\subsection{The subanalytic site}

Let $X$ be a complex analytic curve, denote by $\overline X$ the
complex conjugate curve 
and by $X_\rea$ the underlying real analytic manifold. For $k$ a
commutative ring, we denote by $\mod(k_X)$ the category of sheaves of
$k$-modules on $X$.

We endow $\Op^c_{\xsa}:=\Op^c_{X_{\rea sa}}$ with a Grothendieck
to\-po\-lo\-gy, called the subanalytic topology, by deciding that an usual
open covering $U=\cup_{i\in I}U_i $ in $\Op^c_{\xsa}$ is a covering
for the subanalytic topology if there exists a f\mbox{}inite subset $J\subset
I$ such that $U=\cup_{j\in J}U_j $. Denote by $\xsa$ this site and
call it the \emph{subanalytic site}. Further, denote by $\cov U$ the family of coverings of $U\in\opxsac$ for the subanalytic topology and by
$\mod(k_\xsa)$ the category of sheaves of $k$-modules on the
subanalytic site.

\sloppy One can show (see \cite[Remark 6.3.6]{ast}) that $\mod(k_\xsa)$ is
equivalent to the category of sheaves on the site $\xsalf$, where the class of open sets of $\xsalf$ is $\Op_{\xsa}$ and, for $U\in\Op_\xsa$, the family of coverings of $U$ for $\xsalf$ consists of subanalytic
open coverings $\{U_\sigma\}_{\sigma\in\Sigma}$ of $U$ such that for
any compact $K$ of $X$, there exists a finite subset $J\subset\Sigma$
such that $K\cap\bigr{\cup_{j\in J}U_j}=K\cap U$.

Let $\mathrm{PSh}(k_\xsa)$ be the category of presheaves of
$k$-modules on $\xsa$.
Denote by $for:\mod(k_\xsa)\to \mathrm{PSh}(k_\xsa)$ the forgetful
functor which associates to a sheaf $F$ on $\xsa$ its underlying
presheaf. It is well known that $for$ admits a left adjoint
$\cdot^a:\mathrm{PSh}(k_\xsa)\to \mod(k_\xsa)$.

For $F\in \mathrm{PSh}(k_\xsa)$, let us briefly
recall the construction of $F^a$. 

For $U\in\opxsac$ and  $S=\{U_1,\ldots,U_n\}\in\cov U$, set
\begin{equation}\label{eq_F(S)}
F(S):=\Big\{(s_1,\ldots,s_n)\in\prod_{j=1}^nF(U_j);\ s_j|_{U_j\cap U_k}=s_k|_{U_j\cap U_k},j,k=1,\ldots,n\Big\} \ . \end{equation}
If $S$ is a covering of $U$ and $S'$ is a refinement of $S$, then there exists a natural restriction morphism $F(S)\underset{\rho_{SS'}}{\lra}F(S')$. 

Then, for $U\in\opxsac$, set
\begin{equation}\label{eq_F+} 
F^+(U):=\underset{S\in\cov U}{\varinjlim} F(S) \ .
\end{equation}
It turns out that $F^a\simeq F^{++}$.

\sloppy Now, let $s\in F^a(U)$. Since the inductive limit considered in \eqref{eq_F+} is filtrant, $s$ can be identified to an $n$-uple $(s_1,\ldots,s_n)\in  F(S)$, for $S=\{U_j\}_{j=1,\ldots,n}\in\cov U$, $s_j\in F(U_j)$ ($j=1,\ldots,n$).

Further, if $s\in F^a(U)$ can be identified to $s_1\in F(S_1)$ and to $s_2\in F(S_2)$, for $S_1,S_2\in\cov U$, then there exists a refinement $S\in\cov U$ of $S_1$ and $S_2$ and $\bar s\in F(S)$ such that $s$ can be identified to $\bar s$.

For Proposition \ref{prop_epi_on_xsa} below, see \cite[Proposition
2.1.12]{ast}.

\begin{prop}\label{prop_epi_on_xsa}
Consider the complex in $\mod(k_\xsa)$
\begin{equation}\label{eq_f'ff''}
 F'\overset{\phi}{\longrightarrow} F\overset{\psi}{\longrightarrow}
F'' \ .
\end{equation}
The following conditions are equivalent.
\begin{enumerate}
\item \eqref{eq_f'ff''} is exact.
\item For any $U\in\Op^c_{\xsa}$ and any $t\in F(U)$ such that $\psi(t)=0$,
  there exist $\{U_j\}_{j\in J}\in\cov U$ and $s_j\in F(U_j)$ such that
  $\phi(s_j)=t|_{U_j}$ ($j\in J$).
\end{enumerate}
\end{prop}

We shall denote by
$$ \rho:X\lra\xsa \ ,$$ 
the natural morphism of sites associated to
$\Op^c_{\xsa}\lra\Op_X$. We refer to \cite{ast} for the
definitions of the functors $\rho_*:\mod(k_X)\lra\mod(k_\xsa)$
and $\rho^{-1}:\mod(k_\xsa)\lra\mod(k_X)$ and for Proposition
\ref{prop_functors} below.
\begin{prop}\label{prop_functors}
\begin{enumerate}
\item $\rho^{-1}$ is left adjoint to $\rho_*$.
\item $\rho^{-1}$ has a right adjoint denoted by $\rho_!:\mod(k_X)\lra\mod(k_\xsa)$ .
\item $\rho^{-1}$ and $\rho_!$ are exact, $\rho_*$ is exact on constructible sheaves.
\item $\rho_*$ and $\rho_!$ are fully faithful.
\end{enumerate}
In particular, we can consider $\mod(k_X)$ as a subcategory of $\mod(k_\xsa)$.
\end{prop}

The functor $\rho_!$ is described as follows. If
$U\in\Op^c_\xsa$ and 
$F$ is a sheaf on $X$, then $\rho_!(F)$ is the sheaf on $\xsa$
associated to the presheaf $U\mapsto F\bigr{\overline U}$.


\subsection{Def\mbox{}inition and main properties of $\ot_\xsa$}

Denote by $\D_X$ the sheaf of differential operators with holomorphic
coefficients on $X$. Denote by $\dbxr$ the sheaf of distributions on $X$ and, for a
closed subset $Z$ of $X$, by $\Gamma_Z(\dbxr)$ the subsheaf of
sections supported by $Z$. One denotes by $\dbtxsa$ the presheaf of
\emph{tempered distributions} on $X_{\rea}$ def\mbox{}ined as follows
$$\Op_\xsa\ni U \longmapsto \dbtxsa(U):=\Gamma(X;\dbxr)\big/\Gamma_{X\setminus U}(X;\dbxr) \ .$$
In \cite{ast} it is proved that $\dbtxsa$ is a
sheaf on $\xsa$. This sheaf is well def\mbox{}ined in the category
$\mod(\rho_! \D_X)$. 
Moreover, for any $U\in\opxsac$, $\dbtxsa$ is $\Gamma(U,\cdot)$-acyclic.

One def\mbox{}ines the sheaf $\ot_\xsa \in D^b\bigr{\rho_!\D_X}$ 
of tempered holomorphic functions as
$$\ot_\xsa:=R\mathcal Hom_{\rho_! \mathcal{D}_{\overline X}}\bigr{\rho_!\mathcal{O}_{\overline X},\dbt_{X_{\rea}}}\ .$$

In \cite{ast} it is proved that, since $\dim X=1$, $R\rho_*\O_X$
and $\ot_\xsa$ are concentrated in degree $0$ .
Hence we can write the following exact sequence of sheaves on $\xsa$

$$ 0\longrightarrow \ot_\xsa \longrightarrow \dbtxsa
\overset{\bar\partial}{\longrightarrow} \dbtxsa\longrightarrow 0\ .$$

\begin{lem}
Let $X=\com,\ X_\rea=\rea^2$, $U,V\in\Op^c_{\rea^2_{sa}}$. 
\begin{enumerate}
\item $H^j(U,\ot_\xsa)=0$, for $j>0$.
\item The following sequence is exact
\begin{equation}\label{eq_loj_ot} 0\longrightarrow \ot_\xsa(U\cup V)
  \longrightarrow \qquad\qquad\qquad\qquad\qquad\qquad\qquad\qquad\phantom{a}\end{equation}
$$ \lra \ot_\xsa(U)\oplus\ot_\xsa(V) \longrightarrow\qquad \phantom{a}$$ 
$$\phantom{a}\qquad\qquad\qquad\qquad\qquad\qquad\qquad\qquad\qquad \lra\ot_\xsa(U\cap V)\longrightarrow 0 \ .$$
\end{enumerate}
\end{lem}
\emph{Proof.} \emph{(i)} By the def\mbox{}inition of $\dbtxsa$, 
given $h\in\dbtxsa(U)$, there exists $\widetilde{h}\in\dbxr(\rea^2)$ such that $\widetilde{h}\big|_U=h$. It is well known that there exists $g\in\dbxr(\rea^2)$ such that $\bar{\partial}g=\widetilde h$. This implies that $\bar{\partial}\bigr{g|_U}=h$. So we have the exact sequence
$$ 0\longrightarrow \ot_\xsa(U) \longrightarrow \dbtxsa(U)
\overset{\bar\partial}{\longrightarrow} \dbtxsa(U)\longrightarrow 0 \
. $$

Since $\dbtxsa$ is acyclic with respect to the functor
$\Gamma(U;\cdot)$, for $U\in\Op^c_\xsa$, it follows that, for all $j\in\integergz$,
$H^j(U,\ot_\xsa)=0$.

\emph{(ii)} Obvious from \emph{(i)}.
\proofend

Now we recall the def\mbox{}inition of polynomial growth for
$\mathcal{C}^\infty$ functions on $X_\rea$ and in
\eqref{eq_debar_cinftyt} we give an alternative expression for tempered
holomorphic functions on $U\in\Op_\rsa^c$.

\begin{df}
Let $U$ be an open subset of $X_\rea$,
$f\in\mathcal{C}^\infty_{X_\rea}(U)$. One says that $f$ has
\emph{polynomial growth at $p\in X$} if it satisf\mbox{}ies the
following condition. For a local coordinate system $(x_1,\ldots,x_n)$
around $p$, there exist a suf\mbox{}f\mbox{}iciently small compact
neighborhood $K$ of $p$ and $M\in\integergz$ such that 
\begin{equation}\label{eq_temp}
 \underset{x\in K\cap U}{\sup}\dist x{K\setminus U}^M\big|f(x)\big|<+\infty \ .
\end{equation}
We say that $f\in\mathcal{C}^\infty_{X_\rea}(U)$ has {\em polynomial
  growth on $U$} if it has
polynomial growth at any $p\in X$. We say that $f$ is {\em tempered at
  $p$} if all its derivatives have polynomial growth at $p\in X$. We say
that $f$ is {\em tempered on $U$} if it is tempered at any $p\in X$. Denote
  by $\mathcal C^{\infty,t}_X$ the presheaf on $X_\rea$ of tempered
  $\mathcal C^{\infty}$-functions.
\end{df}

It is obvious that $f$ has polynomial growth at any point of $U$. If
no confusion is possible we will write ``$f$ is tempered" instead of ``$f$
is tempered on $U$".

In \cite{ast} it is proved that $\mathcal C^{\infty,t}_X$ is a sheaf
on $\xsa$.

For $U\subset\rea^n$ a relatively compact open set, we can characterize functions with polynomial growth on $U$ by means of a family of norms. 

For $x\in\rea^n$, $f\in\mathcal{C}^\infty_{\rea^n}(U)$,
$g=(g_1,\ldots,g_m)\in\bigr{\mathcal{C}^\infty_{\rea^n}(U)}^m$ and $M\in\integergz$, denote
\begin{eqnarray}\label{eq_dpart}
\de {\partial U}x & := & \dist x{\partial U} \ ,\\
\norm fMU & := & \underset{x\in U}{\sup}\ \ \delta_{\partial U}(x)^{M} |f(x)| \ ,\notag\\
\norm gMU & := & \max\bigc{\norm {g_j}MU;\ j=1,\dots,m}\notag\ .
\end{eqnarray}

\begin{prop}\label{prop_equiv_def_pol_gr}
Let $U\subset\rea^n$ be a relatively compact open set and let $f\in\mathcal{C}^\infty_{\rea^n}(U)$.
Then $f$ has polynomial growth if and only if there exists $M\in\reagz$ such that
\begin{equation} \label{eq_norm}  \norm fMU<+\infty \ ,\end{equation}
or equivalently: there exist $C,M\in\reagz$ such that for any $x\in U$, 
$$\big|f(x)\big|\leq C\de{\partial U}x^{-M} \ .$$
\end{prop}

\emph{Proof}. 
Suppose that $f$ satisfies \eqref{eq_norm}, that is, 
$$ \underset{x\in U}{\sup}\ \ \dpart U(x)^{M} |f(x)|<+\infty \ .$$ 
Let $K$ be a compact
neighborhood of $\overline U$. For any $p\in\overline U$, $K$ is a compact neighborhood of $p$ such that 
\begin{eqnarray*}
 \underset{x\in K\cap U}{\sup}\dist x{K\setminus U}^M|f(x)|  &  =  &
  \underset{x\in U}{\sup}\ \dpart U(x)^M|f(x)| \\
  &<  & +\infty \ .
\end{eqnarray*}
Hence, $f$ has polynomial growth. 

Conversely, suppose that $h$ has polynomial growth. That is, for
$p\in\partial U$, there exists a compact neighborhood
$K_p$ of $p$ verifying \eqref{eq_temp}. 

Set
$$ V:=\Big\{x\in K_p;\ \delta_{\partial U\setminus K_p}(x)>\dpart U(x)
\Big\}\ . $$
Then for any $x\in V$, $\dpart U(x)=\delta_{\partial U\cap K_p}(x)$.

Since $p\in V$, there exists $\ep\in\reagz$ such that $\overline{\ball p\ep}\subset
V$. Set 
$$Z_p:=\overline{\ball p\ep}\cup(K_p\cap\partial U) \ . $$ 
Then
$Z_p\cap \partial U=K_p\cap\partial U$ and, for any $x\in Z_p\cap U$,
$\dpart U(x)=\delta_{\partial U\cap K_p}(x)=\delta_{\partial U\cap Z_p}(x)$.

Hence,
\begin{eqnarray*}
\underset{x\in Z_p\cap U}{\sup}\dpart U(x)^M |f(x)| &  =  &  \underset{x\in
  Z_p\cap U}{\sup}\delta_{\partial U\cap Z_p}(x)^M |f(x)| \\
  &  \leq  &  \underset{x\in
  K_p\cap U}{\sup}\delta_{\partial U\cap K_p}(x)^M |f(x)| \\
  &  <  &  +\infty \ .
\end{eqnarray*}

Since $\partial U$ is compact, the conclusion follows.
\proofend

Lemma \ref{siu} below is an easy consequence of Cauchy's Formula. See \cite[Lemma 3]{siu}.
\begin{lem}\label{siu}
Let $U$ be a relatively compact open subset of $X$,
$f\in\O_X(U)$ with polynomial growth on $U$. Then $f\in\ot_\xsa(U)$.
\end{lem}

For Proposition \ref{otct} below, see \cite{ast}. 
\begin{prop}\label{otct}
One has the following isomorphism
$$ \ot_\xsa\simeq R\mathcal Hom_{\rho_!\D_{\overline X}}\bigr{\rho_!\mathcal{O}_{\overline X},\mathcal{C}^{\infty,t}_{X_\rea}}\ .$$
\end{prop}

Hence, we deduce the short exact sequence
\begin{equation}\label{eq_debar_cinftyt} 0\longrightarrow \ot_\xsa(U) \longrightarrow
  \mathcal{C}^{\infty,t}_{X_\rea}(U)
  \overset{\bar\partial}{\longrightarrow}
  \mathcal{C}^{\infty,t}_{X_\rea}(U)\longrightarrow 0 \ . \end{equation}


\subsection{Pull-back of tempered holomorphic functions}\label{SUBSECTION_PULLBACK_OT}\label{subsection_pullback_ot}


Recall that, for $U$ a relatively compact open subset of $\rea^n$ and $z\in\rea^n$, we set $\dpart U(z):=\dist z{\partial U}$.

\begin{lem}\label{suban=>polygr}
Let $X$ be an open subset of $\rea^n$, $f:X\to\rea^m$ be a
$\mathcal{C}^\infty$-subanalytic map. Let $U\in\Op^c_{X_{sa}}$,
$V\in\Op^c_{\rea^m_{sa}}$ satisfying 
$f(U)=V$ and $f\bigr{\partial U}=\partial V$. Let
$h\in\mathcal{C}_{\rea^m}^\infty(V)$. 

Then $h$ has polynomial growth
on $V$ if and only if $h\!\circ\! f$ has polynomial growth on $U$.
\end{lem}

\emph{Proof.} Consider the subanalytic continuous functions
$\dpart U,\dpart V\circ\!f|_{\overline U}: \overline U\to\rea_{\geq0}$. Since
$f\bigr{\partial U}=\partial V$ and $f(U)=V$, 
$$ \bigr{\dpart V\circ f|_{\overline U}}^{-1}\bigr{\{0\}}   =   \partial U \ .  $$
In particular, 
$$ \bigr{\dpart V\circ f}^{-1}\bigr{\{0\}}  =    \dpart U^{-1}\bigr{\{0\}} \ . $$

By Theorem \ref{thm_loj_ineq}, there exist
$a,b,\alpha,\beta\in\reagz$ such that, for any $x\in\overline U$, 
\begin{equation}\label{eq_strange1}
a\Bigr{\dpart V\circ f|_{\overline U}(x)}^\alpha\leq\dpart U(x) \ ,
\end{equation}
and
\begin{equation}\label{eq_strange2}
b\bigr{\dpart U(x)}^\beta\leq\dpart V\circ f|_{\overline U}(x) \ .
\end{equation}

\emph{(i)} Suppose that $h\!\circ\!f$ has polynomial growth on
$U$, that is, there exist $C,M\in\reagz$ such that, for any $x\in U$,
$$ \big|h\bigr{f(x)}\big|\leq C\bigr{\dpart U(x)}^{-M} \ . $$

By \eqref{eq_strange1}, we obtain
$$ \big|h\bigr{f(x)}\big|\leq Ca^{-M}\bigr{\dpart V\circ f|_{\overline U}(x)}^{-M\alpha} \ .$$
Since $f(U)=V$, it follows that, for any $y\in V$,
$$ \big|h(y)\big|\leq Ca^{-M}\bigr{\dpart V(y)}^{-M\alpha} \ ,$$
that is, $h$ has polynomial growth on $V$.

\emph{(ii)} Suppose that $h$ has polynomial growth on $V$, that is, there exist
$C',M'\in\reagz$ such that, for any $y\in V$,
$$ |h(y)|\leq C'\bigr{\dpart V(y)}^{-M'}  \ . $$
Since $f(U)=V$, we have, for any $x\in U$,
$$ \big|h(f(x))\big|\leq C'\bigr{\dpart V\circ f(x)}^{-M'}  \ . $$
By \eqref{eq_strange2}, we obtain
$$ \big|h(f(x))\big|\leq C'b^{-M'}\bigr{\dpart U(x)}^{-M'\beta} \ ,$$
that is, $h\circ f$ has polynomial growth on $U$.
\proofend

\begin{thm}\label{thm_pullback_ot}
Let $X$ be an open subset of $\com$, $f\in\mathcal O_\com(X)$
. Let $U\in\opxsac$
such that $f|_{\overline U}$ is an injective map. Let
$h\in\O_X\bigr{f(U)}$. Then, $h\in\ot_{\com_{sa}}\bigr{f(U)}$ if and only
if $h\circ f\in\otxsa(U)$. 
\end{thm}

\emph{Proof.} Since $f$ is an open mapping, $f|_U:U\to f(U)$ is a holomorphic isomorphism.

It is sufficient to prove that $f(\partial U)=\partial\bigr{f(U)}$ in order
to apply Lemma \ref{suban=>polygr}.

\emph{(i)} 
$f(\partial U)\subset\partial\bigr{f(U)}$. For $x\in\partial U$, there exists $\{x_n\}_{n\in\nat}\subset U$ such
that $x_n\underset{\scriptscriptstyle{n\to+\infty}}{\longrightarrow}x$. It follows that
$f(x_n)\underset{{\scriptscriptstyle{n\to+\infty}}}{\longrightarrow}f(x)$, hence
$f(x)\in\overline{f(U)}$. Suppose that $f(x)\in f(U)$. Since $f|_{U}$
is an isomorphism onto $f(U)$, there exists $\overline x\in U$ such
that $f(\overline x)=f(x)$, this contradicts the hypothesis that
$f|_{\overline U}$ is injective. It follows that $f(x)\in\partial\bigr{f(U)}$.

\emph{(ii)} 
$f(\partial U)\supset\partial\bigr{f(U)}$. For
$y\in\partial(f(U))$, there exists $\{y_n\}_{n\in\nat}\subset
f(U)$ such that $y_n\underset{\scriptscriptstyle{n\to+\infty}}{\longrightarrow}y$. Set $x_n:=\big(f|_U\big)^{-1}(y_n)$. Then $\{x_n\}_{n\in\nat}\subset U$ is a bounded sequence. Hence there exists a
subsequence converging to $x\in\overline U$. Since $f(x)=y$ and $f|_{U}$ is an isomorphism onto $f(U)$, $x\in\partial U$.
\proofend

\section{Existence theorem}  \label{SEC_EXIST_THM}

Let $X\subset\com$ be an open neighborhood of $0$, $P$ a differential
operator defined on $X$, whose only possible singular point is $0$. In this
section we study the non-homogeneous ordinary
dif\mbox{}ferential system relative to $P$.

In the f\mbox{}irst subsection we recall some classical results on the
holomorphic solutions of $P$. 

In the second subsection we start by recalling an existence theorem
for tempered holomorphic functions on small open sectors. As said in
the introduction such a result is classical and it has been
treated in more general cases by many authors. We recall the version
obtained by N. Honda in \cite{honda}.
Then we state and prove the main result of this
section which states that given $U\in\Op_\xsa^c$, with $0\in\partial
U$, there exist an open neighborhood $W$ of $0$ and $\bigc{U_j}_{j\in
  J}\in\cov{U\cap W}$ such that $P$ is a surjective endomorphism on
$\otxsa(U_j)$ ($j\in J$). 
The proof is based on the decomposition of the germ of $U$ at $0$ in sets
biholomorphic to open sectors (Theorem \ref{thm_cov_suban}) and on
an existence theorem for sets biholomorphic to open sectors. The proof
of this latter result uses a result on the composition of a
biholomorphism and a tempered holomorphic function (Theorem
\ref{thm_pullback_ot}) in order to reduce the problem to open sectors
of small amplitude.

As a corollary we obtain that $P$ is a surjective
endomorphism of the sheaf $\otxsa$.

\subsection{Some classical results}\label{SUBSEC_CLASSICAL_ODE}

Denote by $\mcC_\com^0$ the sheaf of continuous functions on
$\com$. For $R$ a ring, we denote by $gl(m,R)$ (resp. $GL(m,R)$) the
ring of (resp. multiplicative group of invertible) $m\times m$
matrices. In this chapter we are going to consider $z^{1/l}$,
$l\in\integergz$, as a holomorphic function on open sets contained in open
sectors of amplitude smaller than $2\pi$, by chosing the branch of
$z^{1/l}$ which has positive real values on $\reagz\times\{0\}$.

Let $X\subset\com$ be an open disc centered at $0$. Let
\begin{equation}\label{eq_diff_op} 
 P:=z^N\frac d{dz}I_m+A(z) \  ,\end{equation}
where $m\in\integer_{>0}$, $N\in\nat$, $ A\in gl\bigr{m,\O_\com(X)}$ and
$I_m$ is the identity matrix of order $m$. 

Theorem \ref{thm_sol_structure} below is a fundamental result on ordinary
differential systems. A complete proof of Theorem
\ref{thm_sol_structure} below, originally due to Hukuhara and
Turrittin, is given in \cite{wasow}.

\begin{thm}[See \cite{wasow}]\label{thm_sol_structure}
Let $P$ be the dif\mbox{}ferential operator \eqref{eq_diff_op}. There
exist $l\in\integergz$, a diagonal matrix $\Lambda\in gl\big(m,z^{-1/l}\cdot\com[z^{-1/l}]\big)$ and for any
$\theta_0\in\rea$, there exist $\theta_1,\theta_2\in\rea$,
$\theta_1<\theta_0<\theta_2$, $R,K,M\in\reagz$ and
$F_{\gsect{\theta_1}{\theta_2}R}\in
GL\Big(m,\O_\com(\gsect{\theta_1}{\theta_2}R)\cap\mathscr{C}^0_\com\big(\overline{\gsect{\theta_1}{\theta_2}R}\setminus\{0\}\big)\Big)
$, satisfying the following conditions
\begin{enumerate}
\item for any $z\in\gsect{\theta_1}{\theta_2}R$,
\begin{equation}\label{estimate}
  K^{-1}|z|^{M}\leq\big|F_{\gsect{\theta_1}{\theta_2}R}(z)\big|\leq K|z|^{-M} \ , 
\end{equation}
\item the $m$ columns of the matrix
$F_{\gsect{\theta_1}{\theta_2}R}(z)\exp\bigr{\Lambda(z)}$ are
$\com$-linearly independent holomorphic solutions of the system $Pu=0$.
\end{enumerate}
\end{thm}

If no confusion is possible we will write $F(z)$ instead of
$F_{\gsect{\theta_1}{\theta_2}R}(z)$. 

Note that \eqref{estimate} implies that $F,F^{-1}\in GL\big(m,\otxsa(\gsect{\theta_1}{\theta_2}R)\big)$.

\begin{df}
We call the matrix $F(z)\exp\bigr{\Lambda(z)}$, given in Theorem
\ref{thm_sol_structure}, a {\em fundamental holomorphic solution of
  $P$ on $\gsect{\theta_1}{\theta_2}R$}. If $U$ is an open subset of
$\gsect{\theta_1}{\theta_2}R$, we say that $P$ {\em a fundamental
  holomorphic solution on $U$}.
\end{df}

\begin{lem}\label{lem_hom_sol_sector}
Let $U\in\opxsac$, connected and simply connected.
Suppose that $P$ has a fundamental holomorphic solution
$F(z)\exp(\Lambda(z))$ on $U$. Let $g\in\mathcal{O}(U)^m$, $z_1\in U$.

Then, for $\Gamma_{z_1,z}\subset U$ a path from $z_1$ to $z\in U$,
$$
F(z)\exp\bigr{\Lambda(z)}\int_{\Gamma_{z_1,z}}\!\!\!\!\exp\bigr{-\Lambda(\zeta)}F(\zeta)^{-1}\frac{g(\zeta)}{\zeta^N}
d\zeta    $$
is a holomorphic solution of $Pu=g$.
\end{lem}

\emph{Proof.} Obvious.
\proofend

\subsection{Existence theorem for tempered holomorphic functions}

Let $l\in\integer_{>0}$, $p(z)\in z^{-1/l}\cdot\com[z^{-1/l}]$, $S$ an
open sector of amplitude smaller than $2\pi$, $g\in\O_X(S)$. Set
\begin{equation}\label{eq_integral}
I_{p,z_0}(g)(z):=\exp\bigr{p(z)}\int_{\Gamma_{z_0,z}}\exp\bigr{-p(\zeta)}g(\zeta) d\zeta \ ,
\end{equation}
where $z_0\in\overline{S}$ and $\Gamma_{z_0,z}\subset\overline{S}$ is
a path from $z_0$ to $z\in S$.

\begin{prop}[See \cite{honda}, Proposition 2.3]\label{prop_honda}\label{PROP_HONDA}
Let $l\in\integer_{>0}$ and $p(z)\in
z^{-1/l}\cdot\com[z^{-1/l}]$. There exists $\alpha\in\reagz$ such that
for any open sector $S$ of amplitude $\eta\leq\alpha$, there exist
$z_0\in\overline S$ and a path $\Gamma_{z_0,z}\subset\overline{S}$ from $z_0$ to $z\in S$ such that if $g\in\O_X(S)$ satisf\mbox{}ies
$\norm gM{S}<+\infty$, for some $M\in\reagz$, then
$I_{p,z_0}(g)\in\O_X(S)$ and 
$$\norm {I_{p,z_0}(g)}MS< +\infty \ . $$  
\end{prop}

Now we prove an analogue of Proposition \ref{prop_honda}
for sets biholomorphic to an open sector of sufficiently small
amplitude. Then we will use such a result
to prove an existence theorem for $P$ on $U\in\opxsac$,
$0\in\partial U$.

\begin{prop}\label{prop_Tcase_honda}
Let $W\subset\com$ be an open neighborhood of $0$,
$\phi\in\O_\com(W)$, $\phi(0)=0$, $l\in\integergz$, $p\in
z^{-1/l}\com[z^{-1/l}]$.

There exist $r,\eta\in\reagz$ such that for any open sector
$S\subset\subset\ball0r\subset W$ of amplitude smaller than $\eta$,  there exist 
$z_{0}\in\phi\bigr{\overline S}$ and a path
$\Gamma_{z_{0},z}\subset\phi\bigr{\overline S}$ from $z_{0}$ to
$z\in\phi(S)$ such that, for any $g\in\otxsa\bigr{\phi(S)}$,
$$
I_{p,z_{0}}(g)(z)=\exp\bigr{p(z)}\int_{\Gamma_{z_{0},z}}\!\!\!\!\!\!\exp\bigr{-p(\zeta)}g(\zeta)d\zeta\
\in\ \otxsa\bigr{\phi(S)} \ .  $$
\end{prop}

\emph{Proof.} The proof is based on the following sequence of
equivalences which will be made rigorous along the proof.
\begin{eqnarray}\label{eq_last_line}
I_{p,z_0}(g)(z)  &  \in  &  \otxsa(\phi(S))\notag\\
  &  \Updownarrow  & \notag \\
I_{p,z_0}(g)\circ\phi(w)  &  \in  &  \otwsa(S)\notag\\
  &  \Updownarrow  &  \notag\\
I_{\tilde p,w_0}(\tilde g)(w)  &  \in  &  \otwsa(S)
\end{eqnarray}
for some $\tilde p(w)\in w^{-1/l'}\cdot\com[w^{-1/l'}]$,
$l'\in\integergz$, $w_0\in\overline S$ and $\tilde
g\in\otwsa(S)$. We will obtain \eqref{eq_last_line} from Proposition
\ref{prop_honda} by taking the amplitude of $S$ small enough.

There exists $c\in\integer_{>0}$ such that, $\phi(w)=w^c\phi_1(w)$ and
$\phi_1(0)\neq0$, for any $w\in W$. There exist $r,\eta_0\in\reagz$, $\eta_0<2\pi $,
be such that, for any open sector
$S\subset\subset\ball0r\subset W$ of amplitude smaller than $\eta_0$,
$\phi|_{\overline S}$ is injective. For the rest of the proof, a
sector $S$ will be supposed to have amplitude (resp. of
radius) smaller than $\eta_0$ (resp. $r$).

Let 
$$p(z):=\sum_{j=1}^n{a_j\over z^{j/l}} \ , $$ 
for $q\in\integer_{>0}$ and
$a_j\in\com$ ($j=1,\ldots,n$). 

We have
\begin{eqnarray*}
p\bigr{\phi(w)}  
  &  =  &  \sum_ {j=1}^{n}{a_j\over\bigr{w^{c}\phi_1(w)}^{j/l}}  \\
  &  =  &  \sum_{j=1}^{n}a_j {\phi_{2,j}(w)\over w^{cj/l}}  \\
  &  =  &  \sum_{j=1}^{n}a_j \Bigg(\sum_{k=1}^{q_j}{\beta_{j,k}\over
  w^{k/\lambda_j}}+\phi_{3,j}(w)\Bigg) \\
  &  =  &  \sum_{j=1}^{q'}{a'_j\over w^{j/l'}}+\psi_j(w) \ ,
\end{eqnarray*}
 for some $l',\lambda_j,q_j,q'\in\integer_{>0}$, $\beta_{j,k},a'_j\in\com$
 and $\phi_{2,j},\phi_{3,j},\psi_j$ power series in
 $z^{1/l''}$, for some $l''\in\integergz$, converging on $S$ and
 defined on $\overline{S}$.

Set 
$$ \tilde p(w):=\sum_{j=1}^{q'}{a'_j\over w^{j/l'}}\ \in\  w^{-1/l'}\com[w^{-1/l'}] $$
and
$$ h(w):=\exp\biggr{\sum_{j=1}^{q'}\psi_j(w)}\ \in\ \O_\com(S)\cap\mcC_\com^0\bigr{\overline{S}} \ .  $$

It follows that
$$ \exp\big(p(\phi(w))\big)=\exp\big(\tilde p(w)\big) h(w) \ .$$

Consider $\tilde p\in w^{-1/l'}\com[w^{-1/l'}]$. By Proposition
\ref{prop_honda}, there exists $\eta\in\reagz$, such that for
$S$ an open sector of amplitude smaller that $\eta$,
there exist $w_{0}\in\overline S$, a path
$\Gamma_{w_{0},w}\subset\overline S$ from $w_{0}$ to $w$ such
that, for any $\tilde g\in\otwsa(S)$,
\begin{equation}\label{eq_int_tp}
\exp\bigr{\tilde p(w)}\int_{\Gamma_{w_{0},w}}\!\!\!\!\!\!\!\!\exp\bigr{-\tilde
  p(\zeta)}\tilde g(\zeta)d\zeta\ \in\ \otwsa(S) \ .
\end{equation}

Since the multiplication by $h$ and
$h^{-1}$ is a bijection on $\otwsa(S)$, \eqref{eq_int_tp} implies
that, for any $\tilde g\in\otwsa(S)$,

\begin{eqnarray}\label{eq_int_tp_h}
{\lefteqn {\qquad\quad h(w)I_{\tilde p,w_{0}}\bigr{h^{-1}\!\cdot\tilde g} (w)=}}
\\  
&  &  h(w)\exp\bigr{\tilde p(w)}\int_{\Gamma_{w_{0},w}}\!\!\!\!\!\!\!\!\exp\bigr{-\tilde
  p(\zeta)}h(\zeta)^{-1}\tilde g(\zeta)d\zeta\ \in\ \otwsa(S) \ .\notag
\end{eqnarray}

Set $z_{0}:=\phi(w_{0})\in\phi\bigr{\overline S}$ and let $\Gamma_{z_{0},z}:=\phi\bigr{\Gamma_{w_{0},w}}$. Then, 
for any $g\in\otxsa\bigr{\phi(S)}$,

\begin{equation}\label{eq_p_tp}
I_{p,z_{0}}(g)\circ\phi(w)=h(w)I_{\tilde p,w_{0}}\bigr{h^{-1}\!\cdot\!(g\circ\phi)\!\cdot\!\phi'}(w) \ .
\end{equation}

Up to shrinking $\eta$, we can suppose that $\eta<\eta_0$. In
particular $\phi|_{\overline S}$ is injective for any open sector $S$
of amplitude smaller than $\eta$.

Since $(g\circ\phi)\!\cdot\!\phi'\in\otwsa(S)$, 
\eqref{eq_int_tp_h} and \eqref{eq_p_tp} imply that 
$$ I_{p,z_{0}}(g)\circ\phi(w) \in\otwsa(S) \ . $$

Since $\phi|_{\overline S}$ is injective, the conclusion follows by
Theorem \ref{thm_pullback_ot}.
\proofend

Let us now consider the differential operator $P$ given in
\eqref{eq_diff_op}.

\begin{prop}\label{prop_TsatisfiesP}
Let $J$ be a finite set, $W_j\subset\com$ open neighborhoods of $0$,
$\phi_j\in\O_\com(W_j)$, $\phi_j(0)=0$ ($j\in J$). There exist $r,\eta\in\reagz$ such that for any sector
$S\subset\subset\ball0r\subset\cap_{j\in J}W_j$ of amplitude smaller than $\eta$, 
$$ P:\otxsa(\phi_j(S))^m\lra\otxsa(\phi_j(S))^m $$ 
is an epimorphism ($j\in J$).
\end{prop}

\emph{Proof.}
There exists $\eta_0\in\reagz$ such that for any sector
$S\subset\subset\cap_{j\in J}W_j$ of amplitude smaller than $\eta_0$, $P$ has fundamental holomorphic solutions
$F(z)\exp(\Lambda(z))$ on $\phi_j(S)$, for any $j\in J$.

For $k=1,\ldots,m$, let $p_k\in z^{-1/l}\com[z^{-1/l}]$ be the $(k,k)$-entry of
$\Lambda$, for some $l\in\integergz$. 

By Proposition \ref{prop_Tcase_honda}, for any $j\in J,k=1,\ldots,m$, there
exist $r_{j,k},\eta_{j,k}$ such that for any open sector $S\subset\subset\ball0{r_{j,k}}\subset\cap_{j\in J}W_j$ of
amplitude smaller than $\eta_{j,k}$, there exist $z_{0,j,k}\in\phi_j\bigr{\overline S}$ and paths $\Gamma_{z_{0,j,k},z}$ from
$z_{0,j,k}$ to $z\in\phi_j(S)$ such that for any $g_j\in\otxsa(\phi_j(S))$
 $$  \exp\bigr{p_k(z)}\int_{\Gamma_{z_{0,j,k},z}}\!\!\!\!\!\!\!\!\exp\bigr{-p_k(\zeta)}g_j(\zeta)d\zeta\
\in\ \otxsa\bigr{\phi_j(S)}  \ .  $$
\nopagebreak

Set 
$$r:=\min\{r_{j,k};\, j\in J,k=1,\ldots,m\}$$
$$\eta:=\min\{\eta_0,\eta_{j,k};\, j\in J,k=1,\ldots,m\} \ . $$ 
Let $S\subset\subset\ball0r$ be an open sector of amplitude smaller
than $\eta$. Let $\Gamma_j$ be the collection of $m$ paths
$\Gamma_{z_{0,j,k},z}$, then for any $g\in\otxsa(\phi_j(S))^m$
$$ \exp\bigr{\Lambda(z)}\int_{\Gamma_{j}}\exp\bigr{-\Lambda(\zeta)}g(\zeta)d\zeta\
\in\ \otxsa\bigr{\phi_j(S)}^m \  $$
$(j\in J)$.

Since multiplication by $F(z)$ and $\frac{F(z)^{-1}}{z^N}$ are
bijections on $\otxsa\bigr{\phi(V_j)}^m$, we obtain that. for any
$g\in\otxsa(\phi_j(S))^m$,
$$ F(z) \exp\bigr{\Lambda(z)}\int_{\Gamma_{j}}\exp\bigr{-\Lambda(\zeta)}\frac{F(\zeta)^{-1}}{\zeta^N}g(\zeta)d\zeta\
\in\ \otxsa\bigr{\phi_j(S)}^m $$
($j\in J$).

The conclusion follows by Lemma \ref{lem_hom_sol_sector}.
\proofend

\begin{lem}\label{lem_cap}
Let $U,V\in\opxsac$. If $P$ is a surjective endomorphism both on
  $\otxsa(U)^m$ and $\otxsa(V)^m$. Then $P$ is a surjective
  endomorphism on $\otxsa(U\cap V)^m$.
\end{lem}

\emph{Proof.}
The result follows from the exact sequence \eqref{eq_loj_ot}.
\proofend

\begin{thm}\label{thm_P_epi_cov}
Let $U\in\opxsac$ with $0\in\partial U$. There exist an open neighborhood $W\subset X$ of
  $0$ and  $\{U_j\}_{j\in J}\in\cov{U\cap W}$ such that, for any $j\in
  J$,
$$ P:\otxsa(U_j)^m\lra\otxsa(U_j)^m $$
is an epimorphism.
\end{thm}

\emph{Proof.}
By Theorem \ref{thm_cov_suban}, there exist an open neighborhood $W$
of $0$, a f\mbox{}inite set $J$, open sectors $S_{j,k}$,
$\phi_{j,k}\in\O_\com\bigr{\overline{S_{j,k}}}$, such that
$\phi_{j,k}(0)=0$, $\phi_{j,k}|_{\overline{S_{j,k}}}$ is
injective ($j\in J,k=1,2$) and
\begin{equation}\label{eq_cov_suban}
U\cap W=\bigcup_{j\in
  J}\Big(\phi_{j,1}\big(S_{j,1}\big)\cap\phi_{j,2}\big(S_{j,2}\big)\Big) \ .
\end{equation}
Further, by Remark \ref{rem_refined}, we can suppose that the amplitude
and the radius of $S_{j,k}$ are arbitrarly small. In particular, Proposition
\ref{prop_TsatisfiesP} applies and we have that $P$ is an
epimorphism on $\phi_{j,k}\bigr{S_{j,k}}$, for any $j\in J$, $k=1,2$. 

The conclusion follows from \eqref{eq_cov_suban} and Lemma \ref{lem_cap}.
\proofend

The following Corollary is an obvious consequence of Theorem
\ref{thm_P_epi_cov}. In view of Proposition \ref{prop_epi_on_xsa}, it
states that $P$ is an epimorphism of sheaves on $\xsa$.

\begin{cor}\label{cor_P_epi_sh_xsa}
Let $U\in\opxsac$ with $0\in\partial U$. There exist an open
  neighborhood $W\subset X$ of $0$ such that for any
  $g\in\otxsa(U)^m$, there exist $\{U_j\}_{j\in J}\in\cov{U\cap W}$
  and $u_j\in\otxsa(U_j)^m$ satisfying 
$$ \phantom{(j\in J)}\qquad\qquad\qquad  Pu_j=g|_{U_j} \qquad\qquad\qquad (j\in J) \ . $$
\end{cor}

\emph{Proof.} Obvious.
\proofend

\section{Tempered holomorphic solutions}  \label{SEC_DMOD}

In this section we deal with solutions of $\D_X$-modules, for $X$ a
complex analytic curve. 

In the first subsection we recall some classical results about
$\D_X$-modules. First, for a coherent $\D_X$-module $\M$, we define
the complex of holomorphic (resp. tempered holomorphic) solutions of $\M$,
$\sol\,\M$ (resp. $\solt\M$). Then, we recall that, if $\M$ is a regular holonomic
$\D_X$-module, then $\solt\M\simeq\sol\,\M$. Moreover we recall that a holonomic
$\D_X$-module is locally an extension of a
$\D_X$-module supported on a point (hence regular) and a $\D_X$-module
locally isomorphic to a differential operator.

In the second subsection we state the existence theorem in the
framework of $\D$-modules. It asserts that, for a holonomic
$\D_X$-module $\M$, $H^1\bigr{\solt\M}$ is isomorphic to
$H^1\bigr{\sol\M}$. Using the results recalled in the first
subsection, we reduce to the case of a differential operator. Such
case is the object of the third subsection.

In the third subsection we treat the case of a differential
operator. Making use of the language of sheaves on $\xsa$, we give a
more natural setting and statement to the results obtained in Section
\ref{SEC_EXIST_THM}.

In the fourth subsection we prove that $\solt(\M)$ is $\rea$-constructible in
the sense of sheaves on $\xsa$.

Throughout this section, $X$ will be a complex analytic curve.

\subsection{Classical results on $\D$-modules}

For a detailed and comprehensive exposition of $\D_X$-modules we refer to \cite{bjork} and
\cite{kashi}. For an introduction to derived categories and cohomology
of sheaves, we refer to \cite{som}. 

We denote by $\D_X$ the sheaf of differential operators with
holomorphic coefficients on $X$, $\dmod X$ the category of
$\D_X$-modules, $\cdmod X$ the full subcategory of $\dmod X$
consisting of coherent $\D_X$-modules.

For $\M\in\cdmod X$ we denote by $\ch\M$ the characteristic variety of
$\M$. Recall that $\M\in\cdmod X$ is said \emph{holonomic} if
$\dim\ch\M=1$. We denote by $\hdmod X\subset\cdmod X$ the abelian
category of holonomic $\D_X$-modules.

\sloppy We denote by $D^b(\xsa)$ (resp. $D^b(X)$, $D^b(\D_X)$) the
bounded derived category of sheaves of $\com$-vector spaces on $\xsa$
(resp. sheaves of $\com$-vector spaces on $X$, $\DX$-modules). We
denote by $\dbdx[coh]$ (resp. $\dbdx[h]$) the full subcategory of $\dbdx$
consisting of bounded complexes whose cohomology groups are coherent (resp. holonomic). 
For
$\M\in\dbdx[coh]$, set $ \ch\M:= \cup_{j\in\integer}\ch H^j(\M)$. 

Let $T^*X$ be the cotanget bundle on $X$, $\pi_X:T^*X\to X$ the
canonical projection, $T^*_XX$ the zero section of $T^*X$ and
$\dot T^*X:=T^*X\setminus T^*_XX$. 

For $\M\in\dbdx[coh]$, set 
$$ S(\M):=\pi_X\Bigr{\ch\M\cap\dot T^*X} \ . $$

It is well known that, if $\M\in\dbdx[h]$, then $S(\M)$ is a discrete subset of $X$.

\begin{df}
An object $\M\in\dbdx[h]$ is said \emph{regular holonomic} if, for any $x\in X$,
$$ \mathrm{RHom}_{\D_X}(\M,\O_{X,x})\overset{\sim}{\lra} \mathrm{RHom}_{\D_X}(\M,\widehat\O_{X,x})  \
, $$
where $\widehat\O_{X,x}$ is the $\D_{X,x}$-module of formal power
series at $x$. We denote by $\dbdx[rh]$ the full subcategory of
$\dbdx[h]$  of regular holonomic $\D_X$-modules. 
\end{df}



Recall that $\rho:X\to\xsa$ is the natural morphism of sites. For a
coherent $\D_X$-module $\M$, we set for short
\begin{eqnarray*}
\sol\M  & = & R\rho_*\mathrm R\mathcal Hom_{\D_X}\big(\M,\O_X\big) \in D^b\big(\xsa\big)\\
\solt\M & = & \mathrm R\mathcal Hom_{\rho_!\D_X}\big(\rho_!\M,\ot_\xsa\big) \in D^b\big(\xsa\big)\ .
\end{eqnarray*}

For Theorem \ref{thm_regular_case} below, see \cite{rims}. We recall it here with the notation of \cite{arx}. 

\begin{thm}\label{thm_regular_case}
Let $X$ be a complex analytic manifold, $\M\in\dbdx[rh]$. The na\-tu\-ral morphism in $D^b(\xsa)$
$$\solt\M\longrightarrow\sol\M  $$
is an isomorphism.
\end{thm}

For $a\in X$, let $ \Gamma_{[\{a\}]}(\cdot)$ be the tempered support
functor on $a$. For $\M\in\dbdx$, denote
$$ \lM:=\O_X(*a)\otimes_{\O_X}\M \ ,$$
where $\O_X(*a)$ is the $\D_X$-module of meromorphic functions at $a$.

Proposition \ref{prop_lM_P} below follows from Kashiwara's Lemma (see
\cite[Theo\-rem 4.30]{kashi}) and Kashiwara's thesis \cite{kashi70}.
\begin{prop}\label{prop_lM_P}
Let $a\in X$.
\begin{enumerate}
\item  For $\M\in\dbdx[h]$, there exists a distinguished triangle
$$R\Gamma_{[\{a\}]}\M\lra\M\lra \lM\overset{+1}{\lra} \ .$$
\item If $\M\in\dbdx[h]$, then $R\Gamma_{[\{a\}]}\M\in\dbdx[rh]$.
\item If $\M\in\hdmod X$, then there exist an open neighborhood $W$
  of $a$ and $P\in\D_{X}(W)$ such
that 
$$ \lM|_W\simeq\frac{\D_{X}|_W}{\D_{X}|_W\cdot P} \ .$$
\end{enumerate}
\end{prop}

\subsection{Existence theorem for holonomic $\D_X$-modules}\label{SUBSEC_H1}

\begin{thm}\label{thm_ext1}
Let $\M\in\hdmod X$. The natural morphism of sheaves on $\xsa$
\begin{equation}\label{eq_ext1} 
H^1\Bigr{\solt\big(\M\big)}\lra H^1\Bigr{\sol\bigr{\M}}
\end{equation}
is an isomorphism.
\end{thm}

\emph{Proof.}
The problem is local on $\xsa$. Since $S(\M)$ is a discrete set, 
it is sufficient to prove the statement in the case
$S(\M)\subset\{a\}$, for $a\in X$.

Now, using Theorem \ref{thm_regular_case} and Proposition
\ref{prop_lM_P}, it is sufficient to prove the statement for
$\M=\frac{\D_W}{\D_W\cdot P}$, for $W\subset X$ an open neighborhood of
$a$ and $P\in\D_X(W)$.

That is, up to shrinking $X$, we are reduced to prove that the natural morphism of sheaves on $\xsa$
$$ H^1\Bigr{\solt\Bigr{\frac{\D_X}{\D_X\cdot P}}}\lra H^1\Bigr{\sol\Bigr{\frac{\D_X}{\D_X\cdot P}}} $$
is an isomorphism. 


This is the object of Subsection \ref{subsection_P_case} below.

\proofend

\subsection{The case of a single operator}\label{subsection_P_case}

For this subsection, let $X\subset\com$ be an open disc centered at the
origin and
\begin{equation}\label{eq_operator}
P=\sum_{j=0}^ma_j(z)\frac{d^j}{dz^j} \ ,\end{equation}
 for $a_j(z)\in\O_X(X)$ ($j=1,\ldots,m$), $a_m$ not identically zero. 

Set $S(P):=S\bigr{\D_X/\D_X\cdot P}$, then we have
$$ S(P)=\bigc{z\in X;\ a_m(z)=0} \ . $$

Remark that, since $\rho_*$ is exact on constructible sheaves and
$\O_X$ is $\rho_*$-acyclic,
$$ \frac{\rho_*\O_X}{P\rho_*\O_X}\simeq\rho_*\frac{\O_X}{P\O_X} \ . $$

\begin{prop}\label{prop_P_case}
The natural morphism of sheaves on $\xsa$
\begin{equation}\label{eq_P_case} 
\frac{\otxsa}{P\otxsa}\longrightarrow\rho_*\frac{\O_X}{P\O_X} \ ,
\end{equation}
is an isomorphism.
\end{prop}

We need two preliminary lemmas.

\begin{lem}\label{lem_0notinU}
Let $U\in\opxsac$, $S(P)\cap U=\emptyset$. For any $g\in\ot(U)$,
there exist $\{U_j\}_{j\in J}\in\cov U$ and $u_j\in\ot(U_j)$ such that
$Pu_j=g|_{U_j}$.
\end{lem}

\emph{Proof.}
Since the problem is local on $\xsa$ and $S(P)$ is a discrete set, we
can suppose that $S(P)=\{0\}$.

\emph{First case}: $0\notin\partial U$. It follows that $P$ is a regular operator on a
neighborhood of $\overline U$. The result follows immediatly by
Theorem \ref{thm_regular_case}.

\emph{Second case}: $0\in\partial U$. The result follows from Theorem
\ref{cor_P_epi_sh_xsa} and the first case.

\proofend

\begin{lem}\label{lem_ball}
Let $U\subset X$ be an open ball and assume $\partial U\cap S(P)=\emptyset $.
Then, the natural morphism
$$ \frac{\otxsa(U)}{P\otxsa(U)}\overset{\displaystyle\phi_t}{\lra}\frac{\O_X(U)}{P\O_X(U)} $$
is an isomorphism.
\end{lem}

\emph{Proof.} Consider the following commutative diagram 
$$\xymatrix{
\frac{\O_X(\overline U)}{P\O_X(\overline U)}\ar[r]^{\phi_c}\ar[rd]  &  \frac{\O_X(U)}{P\O_X(U)} \\
  &  \frac{\otxsa(U)}{P\otxsa(U)}\ar[u]_{\phi_t}
\ . }$$

The proof consists of two steps:
\begin{enumerate}
\item $\phi_c$ is an isomorphism,
\item $\phi_t$ is injective.
\end{enumerate}

(i) Consider the complex
$$ \cF:=\quad0\lra\O_{X}\overset{P}{\lra}\O_{X}\lra0 \ . $$

Since $\partial U\cap S(P)=\emptyset$, then $R\Gamma_{\partial
  U}\bigr{\cF|_{\overline U}}\simeq 0$. It follows that $R\Gamma\bigr{\overline
  U,\cF}\overset{\sim}{\lra}R\Gamma\bigr{U,\cF}$. In particular, since
  $\O_X$ is $\Gamma\bigr{U,\cdot}$ and $\Gamma\bigr{\overline U,\cdot}$-acyclic,
  it follows that $\phi_c$ is an isomorphism.






(ii) Let $h\in\ker(\phi_t)$, that is, $h\in\otxsa(U)$ and there exists $u\in\O_X(U)$
satisfying $Pu=h$. Let us prove that $u\in\otxsa(U)$.

The problem is local on $\xsa$.

Clearly, $u|_{U_0}\in\otxsa(U_0)$ for any $U_0\in\Op_{U_{sa}}^c$.

So, let $x\in\partial U$, there exists an open neighborhood $W$ of $x$
such that $S(P)\cap\overline{U\cap W}=\emptyset$. In particular, $P$ is
a regular operator on $\overline{U\cap W}$.

By Theorem \ref{thm_regular_case}, the complex
$$ 0\lra\otxsa|_{U\cap W}\overset{P}{\lra}\otxsa|_{U\cap W}\lra0 $$ 
is concentrated in degree $0$. In particular, there exists
$\{V_j\}_{j\in J}\in\cov{U\cap W}$ and $v_j\in\otxsa(V_j)$ such that
$Pv_j=h|_{V_j}$, that is, $P(v_j-u|_{V_j})=0$. Since $S(P)\cap\overline
V_j=\emptyset$, then $v_j-u|_{V_j}$ extends holomorphically up to the
boundary of $V_j$. That is $v_j-u|_{V_j}=w$, for some $w\in\O_X\bigr{\overline
V_j}$. In particular $u|_{V_j}\in\otxsa(V_j)$. The conclusion follows.

\nopagebreak
\proofend

Now we can prove Proposition \ref{prop_P_case}.

\emph{Proof of Proposition \ref{prop_P_case}.} 

Since $S(P)$ is a
discrete set and the statement is local on $\xsa$, we can suppose that
$S(P)\subset\{0\}$.

We are going to prove that, for any $U\in\Op_\xsa^c$, the natural morphism
\begin{equation}\label{eq_otu_ou} \
\frac{\otxsa}{P\otxsa}(U)\overset{\displaystyle\phi}{\longrightarrow}\frac{\O_X}{P\O_X}(U) 
\end{equation}
is an isomorphism.

Consider the presheaves on $\opxsac$ defined by
\begin{eqnarray*}
\opxsac\ni U  & \longmapsto & F^t(U)  :=  \frac{\otxsa(U)}{P\otxsa(U)}
\ , \\
\opxsac\ni U  & \longmapsto &  F(U)   :=  \frac{\O_X(U)}{P\O_X(U)} \ .  \\
\end{eqnarray*}
Recall that, for a presheaf $G$ on $\xsa$, we denote by $G^{a}$ the
associated sheaf on $\xsa$. We have that
$\frac{\otxsa}{P\otxsa}:=F^{t,a}$ and $\rho_*\frac{\O_X}{P\O_X}\simeq F^{a}$.

Suppose that $0\notin U$. Then $F^a(U)\simeq0$ and
Lemma \ref{lem_0notinU} implies that $F^{t,a}(U)\simeq0$.

Suppose now that $0\in U$. First, let us prove that $\phi$ is surjective.

Recall \eqref{eq_F(S)}.

Let $s\in F^a(U)$. Since the inductive limit considered in \eqref{eq_F+} is filtrant, $s$ can be identified to
$(s_0,\ldots,s_n)\in F(S)$, for $S=\{U_0,\ldots,U_n\}\in\cov U$ and
$s_j\in F(U_j)$ ($j=0,\ldots,n$).

Up to take a refinement, we can suppose that $0\in U_0$ is an open ball,
$s_0\in F^t(U_0)$, $0\notin U_k$, $s_k=0$ and $s_0|_{U_0\cap U_k}=0$ as
an element of $F^t(U_0\cap U_k)$ ($k\neq0$).

It follows that $(s_0,0,\ldots,0)$ defines an element of
$F^t(S)$. In particular, it defines a section $s^t\in F^{t,a}(U)$ such
that $\phi(s^t)=s$. Hence $\phi$ is surjective.

Now, let us show that $\phi$ is injective. Let $s^t\in F^{t,a}(U)$ such that $\phi(s^t)=0$.

As before, $s^t$ can be identified with $(s_0^t,\ldots,s_n^t)\in
 F^t(S)$, for $S=\{U_0,\ldots,U_n\}\in\cov U$ and $s^t_j\in F^t(U_j)$
 ($j=0,\ldots,n$).

Up to take a refinement of $S$, we can suppose that $0\in U_0$ is an open
ball, $0\notin U_k$ and $s^t_k=0$ for $k\neq0$.

That is, $s^t$ can be identified to $(s_0^t,0,\ldots,0)\in F^t(S)$.

Now, let $\phi_t:F^t(U_0)\to F(U_0)$. By Lemma
\ref{lem_ball}, $\phi_t$ is
injective. Clearly, $\phi(s^t)=0$ implies that $\phi_t(s_0^t)=0$. Hence
$s^t_0=0$ and $\phi$ is injective.

\proofend

\subsection{$\rea$-constructibility for tempered
  holomorphic solutions}\label{SUBSEC_R-C}

In the study of classical solution sheaves of $\D$-modules, the
notions of micro-support and $\rea$-constructibility play a central
role. We refer to \cite{som} to definitions and classical results. In
\cite{arx}, M. Kashiwara and P. Schapira defined such notions
in the context of sheaves on $\xsa$. Further, they conjectured some
results on tempered holomorphic solutions of holonomic $\D$-modules
involving $\rea$-constructibility corresponding to classical results
on holomorphic solutions.

Proposition \ref{prop_r-c} below follows from the results obtained in
Section \ref{SEC_EXIST_THM}. It proves, in dimension $1$, a conjecture
from \cite{arx} stating that, for $\M$ a holonomic $\D_X$-module,
$\solt(\M)$ is $\rea$-constructible in the sense of sheaves on $\xsa$.

Denote by $\bdc[\rea-c]{\com_X}$ the full triangulated subcategory of
the bounded derived category of $\mod(\com_X)$ consisting of complexes
whose cohomology groups are $\rea$-constructible. In what follows, for
$F\in\bdc{\com_\xsa}$ and $G\in\bdc[\rea-c]{\com_X}$, we set for short
$$ \RH[\com_X]GF:=\rho^{-1}\RH[\com_\xsa]GF\in\bdc{\com_X} \
 $$
and
$$ \Rh[\com_X]GF:=R\Gamma(X,\RH[\com_X]GF) \ . $$

\begin{prop}\label{prop_r-c}
Let $X$ be a complex curve and let $\M\in\bdc[h]{\D_X}$. Then, for
any $G\in\bdc[\rea-c]{\com_X}$, $\RH[\com_X]{G}{\solt(\M)}\in\bdc[\rea-c]{\com_X}$.
\end{prop}

\emph{Proof.} We may suppose that $X\subset\com$ is an
open ball centered at the origin. By d{\'e}vissage we may suppose that
$\M\simeq\frac{\D_X}{\D_X\,P}$, for $P$ a differential operator as in \eqref{eq_operator} such that
$S(P)\subset\{0\}$. Since the triangulated category $\bdc[\rea-c]{\com_X}$
is generated by the objects $\com_U$, for $U\in\Op_\xsa^c$, we may
assume that $G=\com_U$ for such an $U$.

Let $V\in\Op_\xsa^c$ such that $0\notin\overline V$, then Theorem
\ref{thm_regular_case} implies that $\solt(\M)|_V\simeq\sol(\M)|_V$. In particular, 
$$ \RH[\com_X]{\com_U}{\solt(\M)}|_{X\setminus\{0\}}\simeq
\RH[\com_X]{\com_U}{\sol(\M)}|_{X\setminus\{0\}} \ .$$
It follows that $\RH[\com_X]{\com_U}{\solt(\M)}$ is
weakly-$\rea$-constructible on $X$ and $\rea$-constructible on $X\setminus\{0\}$.


Since $\Ho[\com_X]{\com_U}{\solt(\M)}$ is a subsheaf of
$\Ho[\com_X]{\com_U}{\sol(\M)}$ and $\RH[\com_X]{\com_U}{\solt(\M)}$ is
concentrated in degree $0,1$, it remains to prove that the stalk at
$0$ of $\mathcal Ext^1_{\com_X}(\com_U,\solt(\M))$ has finite dimension. 


Since $\RH[\com_X]{\com_U}{\solt(\M)}$ is weakly-$\rea$-constructible, there
exists an open ball $B$ such that 
\begin{equation} \label{eq_fiber_ball}
\mathcal Ext^1_{\com_X}\bigr{\com_U,\solt(\M)}_0\simeq R^1\mathit{Hom}_{\com_X}\bigr{\com_{U\cap
    B},\solt(\M)} \ .
\end{equation}

Recall that $\solt(\M)$ is represented by the complex
$$ 0\lra\otxsa\overset{P}{\lra}\otxsa\lra0 $$
and that $\otxsa$ is $\Gamma(V,\cdot)$-acyclic for
$V\in\Op_\xsa^c$. It follows that the object $\Rh[\com_X]{\com_U}{\solt(\M)}$ is
represented by the complex 
\begin{eqnarray*}
\Gamma\bigr{U,\solt(\M)}  &  :=  &
0\lra\otxsa(U)\overset{P}{\lra}\otxsa(U)\lra0  \ . \\
\end{eqnarray*}

In particular 
$$ R^1\mathit{Hom}_{\com_X}\bigr{\com_U,\solt(\M)}\simeq
H^1\bigr{\Gamma\bigr{U,\solt(\M)}} \ .$$
We conclude the proof by showing that $H^1\bigr{\Gamma\bigr{U,\solt(\M)}}$ has
finite dimension.

First consider the case $0\in U$. By \eqref{eq_fiber_ball}, we can
suppose that $U$ is an open ball. Then, by Lemma \ref{lem_ball}, we
have  
$$ H^1\bigr{\Gamma\bigr{U,\solt(\M)}}\simeq H^1\bigr{\Gamma\bigr{U,\sol(\M)}}  $$
and the conclusion follows.

Suppose now that $0\in\partial U$. 

By Theorem \ref{thm_P_epi_cov} and Lemma \ref{lem_cap}, there exists a
finite covering $\bigc{U_j}_{j\in J}\in\cov U$ such that, for
any $K\subset J$ 
\begin{equation}\label{eq_H1_UK=0} 
H^1\bigr{\Gamma\bigr{U_K,\solt(\M)}}\simeq 0 \ ,
\end{equation}
where we have set for short $U_K:=\cap_{k\in K}U_k$.

Arguing by induction on $n\geq1$, we are going to prove that, for any $n\geq1$ and $K_1,\ldots,K_n\subset J$,
$$ H^1\bigr{\Gamma\bigr{\cup_{h=1}^nU_{K_h},\solt(\M)}} \  $$
has finite dimension. This will conclude the proof.

If $n=1$, the result follows at once from \eqref{eq_H1_UK=0}.

Suppose now that, for any $K'_1,\ldots,K'_{n-1}\subset J$,
\begin{equation}\label{eq_ind_step} 
\dim\ H^1\bigr{\Gamma\bigr{\cup_{h=1}^{n-1}U_{K'_h},\solt(\M)}}<+\infty \ . 
\end{equation}

Consider $K_1,\ldots,K_n\subset J$ and the following distinguished triangle
\nopagebreak[0]
\begin{equation}\label{eq_dt_cupcap} 
\Gamma\bigr{\cup_{h=1}^{n}U_{K_h},\solt(\M)}
\lra\quad\qquad\quad\qquad\qquad\qquad\qquad\qquad\qquad\phantom{a}
\end{equation}\nopagebreak[0]
$$\lra\Gamma\bigr{\cup_{h=1}^{n-1}U_{K_h},\solt(\M)}\oplus\Gamma\bigr{U_{K_{n}},\solt(\M)}
\lra$$\nopagebreak[0]
$$\phantom{a}\qquad\qquad\qquad\qquad\qquad\qquad\qquad\lra
\Gamma\bigr{\cup_{h=1}^{n-1}U_{K_h}\cap U_{K_{n}},\solt(\M)}\overset{+1}{\lra} \ .$$

Clearly $U_{K_h}\cap U_{K_{n}}=U_{K_h\cup K_{n}}$. Then
\eqref{eq_ind_step} implies that the second and the third term of the
distinguished triangle \eqref{eq_dt_cupcap} have finite dimensional
cohomology groups.

The conclusion follows.
\proofend



\addcontentsline{toc}{section}{References}

\vspace{5mm}

\tiny

\nopagebreak[4]
\begin{center}\begin{tabular}{ll}
Giovanni Morando & \emph{or}\\
Universit{\`a} degli Studi di Padova, & Universit{\'e} Pierre et Marie Curie, \\
Dipartimento di Matematica Pura e Applicata, & Institut de
Math{\'e}matiques de Jussieu, \\
via G. Belzoni 7, 35131, Padova, Italy. & 175 rue du Chevaleret, 75013, Paris, France.
\end{tabular}\end{center}
\nopagebreak[4]

\nopagebreak[4]

\hspace{11mm}\textit{e-mail address:} \texttt{ gmorando@math.unipd.it }


\begin{thebibliography}{99}

\bibitem[BM88-1]
{bier_mil_pisa} E. Bierstone, P. D. Milman: ``The local geometry of analytic mappings''. Dottorato di Ricerca in Matematica (Doctorate in Mathematical Research), ETS Editrice, Pisa, 1988.

\bibitem[BM88-2]
{bier_mil} E. Bierstone, P. D. Milman: \emph{Semi-analytic and subanalytic sets}. Publ. Math. I.H.E.S. 67, 1988, 5--42. \texttt{http://www.numdam.org}

\bibitem[Bj93]
{bjork} J.E.Bj{\"o}rk: ``Analytic {${\D}$}-modules and
    applications''. Mathematics and its Applications, 247, Kluwer
    Academic Publishers Group, Dordrecht, 1993.





\bibitem[C00]
{coste} M. Coste: ``An introduction to o-minimal geometry''. Pisa. Istituti Ed. e Poligraf\mbox{}ici Intern., Universit{\`a} di Pisa, Dipartimento di Matematica, 2000. \\\texttt{http://name.math.univ-rennes1.fr/michel.coste/Enseignement.html} 

\bibitem[vdD98]
{dries} L. van den Dries: ``Tame topology and o-minimal
  structures''. London Math. Soc. Lecture Note Series, 248; Cambridge University Press, 1998.

\bibitem[F81]
{forster} O. Forster: ``Lectures on Riemann
      surfaces''. Graduate Texts in Ma\-the\-ma\-tics, 81; Springer-Verlag, 1981.

\bibitem[H92]
{honda} N. Honda: \emph{On the solvability of ordinary dif\mbox{}ferential equations in the space of distributions}. Jour. Fac. Sc. Univ. Tokyo IA, vol. 39, n. 2, 207--232, september 1992. 

\bibitem[K70]
{kashi70} M. Kashiwara, \emph{Algebraic study of systems
  of partial differential equations}. M{\'e}moires de la Soci{\'e}t{\'e}
  Math{\'e}matique de France. Nouvelle S{\'e}rie, 63, 1995.

\bibitem[K79]
{kashi79} M. Kashiwara: \emph{Faisceaux constructibles et syst{\`e}mes holon{\^o}mes
              d'{\'e}quations aux d{\'e}riv{\'e}es partielles lin{\'e}aires {\`a}
              points singuliers r{\'e}guliers}. S{\'e}minaire
              Goulaouic-Schwartz, 1979--1980 (French), Exp. No. 19, 7,
              {\'E}cole Polytech., Palaiseau, 1980.

\bibitem[K84]
{rims} M. Kashiwara: \emph{The Riemann-Hilbert problem
    for holonomic systems}. Publ. Res. Inst. Math. Sci., vol. 20, n.2,
    319--365, 1984.

\bibitem[K03]
{kashi} M. Kashiwara: ``$\D$-modules and microlocal calculus''. Translations of Mathematical Monographs, 217. American Mathematical Society, Providence, RI, 2003.

\bibitem[KS90]
{som} M. Kashiwara, P. Schapira: ``Sheaves on manifolds''. Grundlehren der Math. Wiss. 292 Springer, 1990.

\bibitem[KS96]
{moderate_and_formal} M. Kashiwara, P. Schapira:
              \emph{Moderate and formal cohomology as\-so\-cia\-ted with constructible
              sheaves}. M{\'e}moires de la Soci{\'e}t{\'e} Math{\'e}matique de France.
              Nouvelle S{\'e}rie, 64, 1996.

\bibitem[KS01]
{ast} M. Kashiwara, P. Schapira:
              \emph{Ind-sheaves}. Ast{\'e}risque, 271,
              2001.

\bibitem[KS03]
{arx} M. Kashiwara, P. Schapira: \emph{Microlocal study of Ind-sheaves
  I: micro-support and regularity}. Ast{\'e}risque 284, 143--164
(2003).
\texttt{arXiv:math.AG/0108065}.

\bibitem[\L59]
{loj_studia} S. \L ojasiewicz: \emph{Sur le probl{\`e}me
      de la division}. Studia Math., 18, 1959, 87--136.

\bibitem[\L65]
{loj_semi-an} S. \L ojasiewicz: Ensembles semi-analytiques, IHES 1965.

\bibitem[\L93]
{loj_fourier} S. \L ojasiewicz: \emph{Sur la g{\'e}om{\'e}trie semi- et sous-analytique}. Ann. Inst. Fourier, Grenoble, 43 (1993), 1575-1595. \texttt{http://www.numdam.org}.

\bibitem[M67]
{malgr_ideals} B. Malgrange: ``Ideals of differentiable
    functions''. Tata Institute of Fundamental Research Studies in
    Mathematics, No. 3. Tata Institute of Fundamental Research. Bombay , 1967.


\bibitem[M91]
{malgr_birk} B. Malgrange: \emph{{\'E}quations
    diff{\'e}rentielles {\`a} coefficients polynomiaux}. Progress in
    Mathematics, 96. Birkh{\"a}user Boston Inc. 1991.

\bibitem[R-S89]
{ramis_sibuya} J.-P. Ramis, Y. Sibuya: \emph{Hukuhara domains and fundamental existence and uniqueness
              theorems for asymptotic solutions of {G}evrey
              type}. Asymptotic Analysis, 2, 1989, 1, 39--94.



\bibitem[Sa00]
{sabbah_ast} C. Sabbah: \emph{{\'E}quations diff{\'e}rentielles {\`a} points singuliers
              irr{\'e}guliers et ph{\'e}nom{\`e}ne de {S}tokes en dimension
              2}. Ast{\'e}risque, 263, 2000.

\bibitem[Si70]
{siu} Y. Siu: ``Holomorphic functions of polynomial growth
  on bounded domains''. Duke Math. J., 37, 1970, 77--84.

\bibitem[W65]
{wasow} W. Wasow: ``Asymptotic expansions for ordinary dif\mbox{}ferential equations''. Pure and Applied Mathematics, Vol. XIV, Interscience Publishers John Wiley \& Sons, Inc., New York-London-Sydney, 1965.

\end{thebibliography}
\end{document}